\documentclass[preprint]{elsarticle} % options: preprint or review

\usepackage{graphicx}
\usepackage{booktabs}

\makeatletter
	\def\ps@pprintTitle{%
 	\let\@oddhead\@empty
	\let\@evenhead\@empty
	\def\@oddfoot{\centerline{\thepage}}%
	\let\@evenfoot\@oddfoot}
\makeatother

\usepackage{etoolbox}
\patchcmd{\MaketitleBox}{\footnotesize\itshape\elsaddress\par\vskip36pt}{\footnotesize\itshape\elsaddress\par\parbox[b][36pt]{\linewidth}{\vfill\hfill\textnormal{\today}\hfill\null\vfill}}{}{}%
\patchcmd{\pprintMaketitle}{\footnotesize\itshape\elsaddress\par\vskip36pt}{\footnotesize\itshape\elsaddress\par\parbox[b][36pt]{\linewidth}{\vfill\hfill\textnormal{\today}\hfill\null\vfill}}{}{}%

\usepackage{tikz}
\usepackage{enumitem}
\usepackage{nomencl, multicol}
\usepackage{tikz-imagelabels}
\usepackage{pgfplots} 
\usepackage{pgfgantt}
\usepackage{pdflscape}
\pgfplotsset{compat=newest} 
\pgfplotsset{plot coordinates/math parser=false} 
\newlength\fwidth
\newlength\fheight

\usetikzlibrary{shapes.geometric, arrows, positioning, fit, calc}
\usepackage{amsmath, amssymb}
\usepackage[%pdfmark,%
            colorlinks=true,%
            breaklinks=true,%
            linkcolor=blue,
            urlcolor=blue,%
            citecolor=blue,%
            pdftitle={Robust Design Optimization with Limited Data for a Char Combustion Process}, 
            pdfkeywords={Keywords},	
            pdfauthor={Authors},		
            bookmarksopen=false,
            pdfpagemode=UseNone]{hyperref}
            
\usepackage[margin = 1.2in]{geometry}
\usepackage[T1]{fontenc}
\usepackage[english]{babel}
\usepackage[latin1]{inputenc} % Can cause compiling issues
\usepackage{mathtools}
\usepackage{amsmath}
\usepackage{amsfonts}
\usepackage{amsthm}
\usepackage{amssymb}
\usepackage{array}
\usepackage{algorithm} %ctan.org\pkg\algorithms
\usepackage{algorithmicx}
\usepackage{algpseudocode}
\usepackage{subcaption}
\usepackage{siunitx}

\usepackage{color}
\usepackage{url}
\usepackage{cleveref}
\usepackage{graphicx}
\usepackage{breakcites}
\usepackage{soul} % for using the \ul command with linebreak
\usepackage{enumitem}
\usepackage[title,titletoc,toc]{appendix}

\usepackage{makecell}

\usepackage[textsize=tiny]{todonotes}

{\noindent {\textbf{Proof}:} }%
{\hfill $\Box$ \\[1ex] }

\newcommand{\bit}{\begin{itemize}}
	\newcommand{\eit}{\end{itemize}}
\newcommand{\ben}{\begin{enumerate}}
	\newcommand{\een}{\end{enumerate}}

\usetikzlibrary{shapes.geometric, arrows, shadows, fit}

\graphicspath{{./figures/}}

\makenomenclature

\renewcommand*{\nompreamble}{\begin{multicols}{2}}
\renewcommand*{\nompostamble}{\end{multicols}}
\begin{document}
	
	% Here is a frontmatter of an article I wrote, please just edit
	\begin{frontmatter}
		
		\title{Robust design optimization for a nonlinear system via Bayesian neural network enhanced polynomial dimensional decomposition}
				
		\author[hyu1]{{Hyunho Jang}\corref{first}}
		
		\author[hyu]{{Dongjin Lee}\corref{cor1}}
		
		\ead{dlee46@hanyang.ac.kr; +1 858-246-5327}
		%\cortext[first]{These two authors contributed equally to this work.}
		\cortext[cor1]{Corresponding author}
						
		\address[hyu1]{ Department of Automotive Engineering (Automotive-Computer Convergence), Hanyang University, Seoul, South Korea}
  	    \address[hyu] {Department of Automotive Engineering, Hanyang University, Seoul, South Korea}
		
		\begin{abstract}
            Uncertainties such as manufacturing tolerances cause performance variations in complex engineering systems, making robust design optimization (RDO) essential. However, simulation-based RDO faces high computational cost for statistical moment estimation, and strong nonlinearity limits the accuracy of conventional surrogate models. This study proposes a novel RDO method that integrates Bayesian neural networks (BNN) with polynomial dimensional decomposition (PDD). The method employs uncertainty-based active learning to enhance BNN surrogate accuracy and a multi-point single-step strategy that partitions the design space into dynamically adjusted subregions, within which PDD analytically estimates statistical moments from BNN predictions. Validation through a mathematical benchmark and an electric motor shape optimization demonstrates that the method converges to robust optimal solutions with significantly fewer function evaluations. In the ten-dimensional benchmark, the proposed method achieved a $99.97\%$ mean reduction, while Gaussian process-based and Monte Carlo approaches failed to locate the global optimum. In the motor design problem, the method reduced cogging torque by $94.75\%$ with only 6644 finite element evaluations, confirming its computational efficiency for high-dimensional, strongly nonlinear engineering problems.
		\end{abstract}	
		
		
		% Sometimes needed for SIAM papers, need to choose from AMS website
		%\begin{AMS}
		%         35R60,  % Partial differential equations-  	PDEs with randomness,
		%         62H12,  %Statistics - Multivariate analysis - Estimation
		%         65G99,  % Numerical analysis - Error analysis and interval analysis
		%         65Y20  % Numerical analysis -  Complexity and performance of numerical algorithms
		%\end{AMS}
		
\end{frontmatter}

%%%%%%%%%%%%%%%%%%%%%%%%%%%%%%%%%%%%%%%%
%
\section{Introduction} \label{sec:intro}
%%%%%%%%%%%%%%%%%%%%%%%%%%%%%%%%%%%%%%%%
    Robust design optimization (RDO) derives designs insensitive to performance variations by minimizing the impact of uncertainty on system outputs. Since Taguchi's work~\cite{taguchi_taguchi_1993}, industries including aerospace, automotive, civil, and electronics have widely adopted RDO as a key technology to enhance product reliability~\cite{ mourelatos2006methodology, zaman2011robustness, park2006robust, yao2011review, kim4883196deep}. In practice, engineering design processes face inherent uncertainties---such as manufacturing tolerances and material variations---that degrade the performance of nominal designs. 
    These uncertainties are categorized into two primary types: epistemic uncertainty, which arises from model deficiencies in regions with insufficient training data and can be reduced by acquiring more data; and aleatoric uncertainty, which arises from inherent variability, such as manufacturing tolerances, and remains irreducible.  

    A commonly used RDO formulation includes objective and/or constraint functions using second-order statistical moments, such as the mean and variance of design and uncertain variables. This inherently relies on accurate estimation of these moments. Monte Carlo simulation (MCS) can estimate statistical moments but typically requires a large number of samples to achieve reliable convergence, resulting in high computational cost.   
    To mitigate this, various surrogate models, including polynomial chaos expansion~\cite{lee2020practical}, Kriging~\cite{lee2024novel}, support vector machines~\cite{cho2023performance}, reduced order modeling~\cite{guo2025risk}, artificial neural networks~\cite{kim4883196deep}, and more~\cite{lee2023multifidelity, kim2025interpolation,hyun2026robust}, have been used. However, these models often struggle to maintain predictive accuracy when encountering highly nonlinear responses with high-dimensional inputs. 
    High-dimensional inputs make surrogate model construction computationally prohibitive due to the curse of dimensionality, where the required sample number increases exponentially. These models are generally classified into deterministic and probabilistic approaches. 
    Deterministic surrogate models, such as standard artificial neural networks or polynomial chaos expansions, provide only point estimates and fail to account for epistemic uncertainty or prediction reliability. In contrast, probabilistic models quantify such prediction uncertainty. Gaussian process (GP) models---such as Kriging---offer inherent uncertainty quantification, but they become computationally inefficient for design problems with high-dimensional inputs because their training cost scales cubically with the sample number. Bayesian neural network (BNN) has emerged as a scalable alternative, addressing these computational bottlenecks while maintaining robust  uncertainty quantification.  

    BNN offers a scalable probabilistic modeling approach for design problems with high-dimensional inputs and has been successfully applied across diverse engineering domains~\cite{ngartera2024application, hauzenberger2025bayesian, song2025ferroelectric}. Unlike deterministic neural networks, BNN treats both network weights and biases as probability distributions rather than single-point estimates, thereby enabling explicit quantification of the model's epistemic uncertainty.
    This probabilistic formulation supports active learning strategies that prioritize sampling in high-variance regions, improving data efficiency and reliability at reduced computational cost. By leveraging a deep learning architecture, BNN captures complex variable interactions through hierarchical nonlinear mappings, mitigating the exponential parameter growth typical of high-dimensional problems. That said, maintaining accuracy across an expanding design space remains challenging, as the growing sample requirements still impose a significant computational burden.

    To address the computational cost of iterative moment estimation in RDO, polynomial dimensional decomposition (PDD) surrogate offers an efficient alternative. PDD evaluates statistical moments via closed-form expressions derived from its expansion coefficients, scaling effectively by decomposing responses into low-order component functions~\cite{rahman2008polynomial}. However, this approach struggles with strongly nonlinear responses---such as cogging torque in electric motors---where polynomial bases fail to capture sharp local variations. Moreover, surrogates trained solely within the random variable space lack the extrapolation capability needed to maintain accuracy across the broader design space inherent to RDO \cite{guo2025robust}. A recent work on multi-point single-step (MPSS) RDO~\cite{lee2021robust} addresses this limitation by partitioning the design space into subregions and solving localized RDO subproblems in a divide-and-conquer fashion. That said, MPSS RDO relies on gradient-based optimization, which becomes problematic for responses exhibiting severe nonlinearity---such as cogging torque in electric motors---where reliable gradient information is difficult to obtain, limiting its convergence and solution quality.

    This paper proposes a novel RDO method that integrates BNN with PDD. The combined method leverages BNN as a surrogate to replace expensive model evaluations required for PDD coefficient calculation. This enables rapid estimation of statistical moments at significantly reduced computational cost. 
    The proposed method maximizes sampling efficiency through an active learning strategy driven by BNN predictive uncertainty. The method adopts a MPSS method that partitions the design space into subregions and constructs local surrogate models, while using a global search within each subregion to handle severe nonlinearity. A dynamic resizing strategy further adjusts subregion boundaries during optimization to avoid entrapment in local optima and enhance search robustness. PDD is then applied to the resulting surrogate, exploiting the structural decomposition of component functions for efficient and accurate probabilistic analysis. The key contributions of this work are:
    \begin{itemize} 
    \item[(1)] a novel integration of BNN and PDD that simultaneously achieves computational efficiently and prediction accuracy for nonlinear RDO problems with high-dimensional inputs; 
    \item[(2)] an active learning sampling strategy that leverages predictive uncertainty to progressively enhance surrogate model reliability; 
    \item[(3)] a dynamic subregion resizing strategy combined wiht global search that prevents local optima and promotes robust convergence to the global solution.
    \end{itemize}

    The paper is organized as follows. Section~\ref{section:Theoretical background} reviews the theoretical foundations underlying the proposed method.
    Section~\ref{section:RDO for BNN-PDD} formulates the RDO problems using transformed random variables to maximize computational efficiency of PDD basis functions.
    Section~\ref{section:MP} details the step-by-step procedure of the overall optimization algorithm integrating the described components. Section~\ref{section:Numerical examples} validates the proposed method through a mathematical benchmark and an electric motor design problem. Section~\ref{section:conclusion} concludes the paper. 
% 

%%%%%%%%%%%%%%%%%%%%%%%%%%%%%%%%%%%%%%%%%%%
\section{Theoretical background} \label{section:Theoretical background}
%%%%%%%%%%%%%%%%%%%%%%%%%%%%%%%%%%%%%%%%%%%
%
    This section presents the theoretical foundations of the proposed method. Section~\ref{section:input random variables}  defines the input and output random variables. Section~\ref{section:RDO definition} formulates the RDO problem. Section~\ref{section:PDD} summarizes PDD and
    details the calculation of statistical moments via its expansion coefficients. Section~\ref{section:BNN} describes BNN and its uncertainty quantification ability.
%

%%%%%%%%%%%%%%%%%%%%%%%%%%%%%%%%%%%%%%%%%%%
%
\subsection{Input and output random variables} \label{section:input random variables}
%%%%%%%%%%%%%%%%%%%%%%%%%%%%%%%%%%%%%%%%%%%
    Consider a probability space $(\Omega, \mathcal{F}, \mathbb{P})$, where $\Omega$ denotes the sample space, $\mathcal{F}$ represents a $\sigma$-algebra on $\Omega$, and $\mathbb{P}: \mathcal{F} \rightarrow[0,1]$ is a probability measure. We define an $N$-dimensional input random vector $\mathbf{X} := (X_1, \ldots, X_N)^{\intercal}$ as a mapping from $(\Omega, \mathcal{F})$ to a measurable space $(\mathbb{A}^{N}, \mathcal{B}^{N})$. Here, $\mathbb{A}^{N} \subseteq \mathbb{R}^{N}$ denotes a subset of the real space, and $\mathcal{B}^{N}$ denotes the Borel $\sigma$-algebra on $\mathbb{A}^{N}$.
    
    We assume the input random variables to be mutually independent. We define the joint cumulative distribution function (CDF) of $\mathbf{X}$ as $F_{\mathbf{X}}(\mathbf{x}):=\mathbb{P}[\cap_{i=1}^{N}\{ X_{i} \leq x_{i} \}]$, and the corresponding joint probability density function (PDF) as $f_{\mathbf{X}}(\mathbf{x}):=\frac{\partial^{N}F_{\mathbf{X}}(\mathbf{x})}{\partial x_{1} \cdots \partial x_{N}}$. Under this independence assumptions, we factorize the joint PDF as $f_{\mathbf{X}}(\mathbf{x}):=\prod_{i=1}^{N}f_{i}(X_i)$, where $f_{i}(X_i)$ denotes the marginal PDF of the $i$-th random variable $X_i$. We then express the probability measure induced by $\mathbf{X}$ as $(\mathbb{A}^{N}, \mathcal{B}^{N}, f_{\mathbf{X}}(\mathbf{x})\mathrm{d}\mathbf{x})$.
%

%%%%%%%%%%%%%%%%%%%%%%%%%%%%%%%%%%%%%%%%%%%
%
%\subsection{Output random variables} \label{section:output random variables}
%%%%%%%%%%%%%%%%%%%%%%%%%%%%%%%%%%%%%%%%%%%
    Given the input random vector $\mathbf{X}$ on $(\Omega, \mathcal{F}, \mathbb{P})$ and a model function $y : \mathbb{A}^{N} \rightarrow \mathbb{R}$, we define the output random variable as $Y = y(\mathbf{X}) : (\Omega, \mathcal{F}) \rightarrow (\mathbb{R}, \mathcal{B})$. We assume $Y$ has a finite second moment, ensuring that both the mean and variance of $Y$ are well-defined. The output $Y$ captures the propagation of input uncertainty through the model $y(\cdot)$.
%

%%%%%%%%%%%%%%%%%%%%%%%%%%%%%%%%%%%%%%%%%%%
%
\subsection{Robust design optimization} \label{section:RDO definition}
%%%%%%%%%%%%%%%%%%%%%%%%%%%%%%%%%%%%%%%%%%%
    RDO seeks optimal designs by simultaneously accounting for the mean and variance of performance functions under aleatoric uncertainty. We define the design vector $\mathbf{d} = (d_1,\ldots,d_M)^{\intercal} \in \mathcal{D} \subseteq \mathbb{R}^M$ as the mean values of a selected subset of input random variables $\mathbf{X}$, where $M \leq N$. Here, $\mathbf{X}$ represents aleatoric uncertainty arising from inherent variability in the system. This definition reflects the fact that the mean corresponds to the nominal value---a controllable quantity in actual manufacturing. We generally treat the standard deviations as  fixed parameters determined by manufacturing capability. However, in specific situations such as tolerance design, we can treat standard deviations as additional design variables.

    The PDF of the input random variables $\mathbf{X}$ depends on the design vector $\mathbf{d}$, which we express as $f_{\mathbf{X}}(\mathbf{x}; \mathbf{d})$. We formulate the general RDO problem as 
    \begin{align*}
        \min\limits_{\mathbf{d} \in \mathcal{D}\subseteq\mathbb{R}^{M}} \quad c_{0}(\mathbf{d}) :=& g_{0} (\mathbb{E}_{\mathbf{d}}[y_{0}(\mathbf{X})], \mathbb{V}\mathrm{ar}_{\mathbf{d}}[y_{0}(\mathbf{X})]),\\
        \text{subject to} \quad c_{l}(\mathbf{d}) :=& g_{l}(\mathbb{E}_{\mathbf{d}}[y_{l}(\mathbf{X})], \mathbb{V}\mathrm{ar}_{\mathbf{d}}[y_{l}(\mathbf{X})]) \leq 0, \\ 
        \qquad \qquad l =& 1,\ldots, K, \\
        \quad d_{k,L} \leq &d_{k} \leq d_{k,U}, \quad k=1,\ldots, M,
    \end{align*} 
    where $\mathbb{E}_{\mathbf{d}}[y_{l}(\mathbf{X})] := \int_{\mathbb{R}^N}y_{l}(\mathbf{x})f_{\mathbf{X}}(\mathbf{x};\mathbf{d}) \mathrm{d} \mathbf{x}$ denotes the mean of $y_{l}(\mathbf{X})$, and $\mathbb{V}\mathrm{ar}_{\mathbf{d}}[y_{l}(\mathbf{X})] := \mathbb{E}_{\mathbf{d}}[(y_{l}(\mathbf{X}) - \mathbb{E}_{\mathbf{d}}[y_{l}(\mathbf{X})])^2]$ represents the variance of $y_{l}(\mathbf{X})$. The functions $g_0(\cdot)$ and $g_l(\cdot)$ represent arbitrary compositions of the mean and variance, and $d_{k,L}$ and $d_{k,U}$ denote the lower and upper bounds of each design variable $d_k$, respectively.
    
    In practice, we specify this general RDO formulation into concrete mathematical forms.    We adopt the weighted sum method, a widely approach, formulated as
    \begin{align}
        \min\limits_{\mathbf{d} \in \mathcal{D}\subseteq\mathbb{R}^{M}} \quad c_{0}(\mathbf{d}) := w_{1} \frac{\mathbb{E}_{\mathbf{d}}[y_{0}(\mathbf{X})]}{\mu_{0}} + w_{2} \frac{\sqrt{\mathbb{V}\mathrm{ar}_{\mathbf{d}}[y_{0}(\mathbf{X})]}}{\sigma_{0}},
        \label{eq:weighted sum}
    \end{align}
    where $w_1, w_2 \geq 0$ are weights satisfying $w_1 + w_2 =1$, and $\mu_{0}$ and $\sigma_{0}$ are non-zero, real-valued scaling factors. We typically express the constraints as
    \begin{align}
        c_{l}(\mathbf{d}) := \alpha_l \sqrt{\mathbb{V}\mathrm{ar}_{\mathbf{d}}[y_{l}(\mathbf{X})]} - \mathbb{E}_{\mathbf{d}}[y_{l}(\mathbf{X})] \leq 0, \quad l=1, \ldots, K,
        \label{eq:constraints}
    \end{align}
    where $\alpha_l$ are non-negative, real-valued constants associated with the probabilities of constraint satisfaction. Each constraint enforces robustness by guaranteeing that the system performance remains within the feasible region despite inherent variability.
%
%%%%%%%%%%%%%%%%%%%%%%%%%%%%%%%%%%%%%%%%%%%
%
\subsection{Polynomial dimensional decomposition} \label{section:PDD}
%%%%%%%%%%%%%%%%%%%%%%%%%%%%%%%%%%%%%%%%%%%
    PDD uses the analysis-of-variance (ANOVA) decomposition to analyze statistical moments of multivariate functions. The ANOVA decomposition expresses a multivariate function as a sum of component functions of increasing interaction order---comprising a constant term, univariate effects, and multivariate interactions among input variable subsets. PDD can represents any square-integrable function on the probability space $(\Omega, \mathcal{F}, \mathbb{P})$ through an infinite Fourier-like polynomial expansion~\cite{rahman2008polynomial}. For computational efficiency, we adopt the $S$-variate, $m$-th order PDD approximation. This truncates the expansion by retaining only interaction orders up to $S$ and polynomial orders up to $m$. We define the truncated PDD approximation  $\tilde{y}_{S,m}(\mathbf{X})$ as 
    \begin{align}
    \tilde{y}_{S,m}(\mathbf{X}) := y_{0} + \sum\limits_{\substack{\emptyset \neq u \subseteq \{1,\ldots,N \} \\ 1 \leq |u |\leq S}} \sum\limits_{\substack{\mathbf{j}_{|u|} \in \mathbb{N}^{|u|} \\ |u| \leq |\mathbf{j}_{|u|}| \leq m}} c_{u\mathbf{j}_{|u|}} \psi_{u\mathbf{j}_{|u|}}(\mathbf{X}_{u}) \approx y(\mathbf{X}),
    \label{eq:truncated_PDD}
    \end{align}
    where $y_{0} = \int_{\mathbb{A}^{N}} y(\mathbf{X}) f_{\mathbf{X}}(\mathbf{x}) \mathrm{d} \mathbf{x}$ denotes the mean of the function. In~\eqref{eq:truncated_PDD}, $u \subseteq \{ 1,\ldots,N \}$ denotes a non-empty subset of the variable indices with cardinality $|u|$, and $\mathbf{j}_{|u|} = (j_{i_1}, \ldots, j_{i_{|u|}}) \in \mathbb{N}^{|u|}$ is a multi index whose component $j_{i_k}$ specifies the polynomial order for variable $X_{i_k}$. We define the expansion coefficient $c_{u \mathbf{j}_{|u|}}$ and the multivariate basis function $\psi_{u \mathbf{j}_{|u|}}(\mathbf{X}_{u})$ as 
    \begin{align*}
        c_{u \mathbf{j}_{|u|}} :=& \int_{\mathbb{A}^{N}} y(\mathbf{x}) \psi_{u \mathbf{j}_{|u |}}(\mathbf{x}_{u}) f_{\mathbf{X}}(\mathbf{x})\mathrm{d}\mathbf{x}, \\
        \psi_{u \mathbf{j}_{|u|}} (\mathbf{X}_{u}) =& \prod\limits_{i \in u} \psi_{i,j_{i}} (\mathbf{X}_{i}).
    \end{align*}
    Here, $\psi_{i,j_{i}}$ is a univariate polynomial of order $j_i$ that satisfies orthonormality with respect to the marginal probability density function $f_{X_i}(x_i) \mathrm{d} x_i$ of variable $X_i$. These basis functions satisfy the following orthonormality condition:
    \begin{align*}
        \mathbb{E}[\psi_{u \mathbf{j}_{|u|}}(\mathbf{X}_{u})] &= 
        \begin{cases}
        1, \quad \mathbf{j}_{|u|} = \mathbf{0}, \\
        0, \quad \text{otherwise}.
        \end{cases}\\
        \mathbb{E}[\psi_{u \mathbf{j}_{|u|}}(\mathbf{X}_{u})\psi_{v \mathbf{k}_{|v|}}(\mathbf{X}_{v})] &=  \begin{cases}
        1, \quad u=v \text{ and } \mathbf{j}_{|u|} = \mathbf{k}_{|v|}, \\
        0, \quad \text{otherwise},
        \end{cases}
    \end{align*}
    where $v\subseteq\{1,\ldots,N\}$ is another non-empty set of variable indices and $\mathbf{k}_{|v|}\in\mathbb{N}^{|v|}$ is the corresponding multi-index, analogous to $u$ and $\mathbf{j}_{|u|}$.
    
    For notational convenience, we re-index the multi-index basis functions into a single index $\{ \psi_{i}(\mathbf{X}) \}_{i=1}^{L_{N,S,m}}$. The total number of basis functions $L_{N,S,m} = 1 + \sum\limits_{s=1}^{S} \binom{N}{s} \binom{m}{s}$ depends on the truncation parameters $S$ and $m$. We then express the truncated PDD as a linear combination:
    \begin{align*}
        \tilde{y}_{S,m}(\mathbf{X}) = \sum\limits_{i=1}^{L_{N,S,m}} c_i \Psi_i (\mathbf{X}),
    \end{align*}
    with the expansion coefficient vector $\mathbf{c} = (c_1, \ldots, c_{L_{N,S,m}})^{\intercal} \in \mathbb{R}^{L_{N,S,m}}$. We estimate $\mathbf{c}$ via least squares by evaluating the function at $\bar{L}$ sample points, yielding the linear system: 
    \begin{align*}
        \underbrace{
        \begin{bmatrix}
        \psi_1(\mathbf{x}^{(1)}) & \cdots & \psi_{L_{N,S,m}}(\mathbf{x}^{(1)}) \\ \vdots & \ddots & \vdots \\ \psi_{1}(\mathbf{x}^{(\bar{L})}) & \cdots & \psi_{L_{N,S,m}}(\mathbf{x}^{(\bar{L})})
        \end{bmatrix}}_{=: \mathbf{A} \in \mathbb{R}^{\bar{L} \times L_{N,S,m}}
        }
        \underbrace{
        \begin{bmatrix}
        c_1 \\ \vdots \\c_{L_{N,S,m}}
        \end{bmatrix}}_{=: \mathbf{c} \in \mathbb{R}^{L_{N,S,m}}
        }
        = \underbrace{
        \begin{bmatrix}
        y(\mathbf{x}^{(1)}) \\ \vdots \\ y(\mathbf{x}^{(\bar{L})})
        \end{bmatrix}}_{=: \mathbf{b} \in \mathbb{R}^{\bar{L}}},
    \end{align*}
    from which we obtain 
    \begin{align*}
        \mathbf{c} = (\mathbf{A}^{\intercal} \mathbf{A})^{-1} \mathbf{A}^{\intercal} \mathbf{b}.
    \end{align*}
    Exploiting the orthonormality of the basis functions, we analytically derive the mean and variance~\cite{rahman2010statistical} of $\tilde{y}_{S,m}(\mathbf{X})$ directly from the expansion coefficients:
    \begin{align}
        \mathbb{E} [\tilde{y}_{S,m}(\mathbf{X})] = c_{1}, \quad \mathbb{V}\mathrm{ar} [\tilde{y}_{S,m}(\mathbf{X})] = \sum\limits_{i=2}^{L_{N,S,m}} c^{2}_{i}.
        \label{eq:pdd_mean_var}
    \end{align}
%

%%%%%%%%%%%%%%%%%%%%%%%%%%%%%%%%%%%%%%%%%%%
%
\subsection{Bayesian neural networks}\label{section:BNN}
%%%%%%%%%%%%%%%%%%%%%%%%%%%%%%%%%%%%%%%%%%%
%
    Artificial neural networks model complex nonlinear relationship between an input vector $\mathbf{x}$ and output $y(\mathbf{x})$ through multi-layer neuron structures, i.e., $y(\mathbf{x})\simeq \tilde{y}_{\cal{N}}(\mathbf{x}; \boldsymbol{\theta})$. Here, $\boldsymbol{\theta}=(\omega_1,\ldots,\omega_{L_{\omega}},b_1,\ldots,b_{L_b})^{\intercal}$ denotes the model parameters. Also, $\omega_i$ for $i = 1,\ldots, L_{\omega}$ are the weights and $b_j$ for $j = 1,\ldots, L_b$ are the biases, both determined through training. A conventional deterministic neural network identifies a single optimal parameter set that minimizes a loss function, such as the mean squared error between predicted and observed values. This point estimation approach fails to capture the epistemic uncertainty inherent in the model itself.
    \begin{figure}[htb!]
        \centering
        \begin{minipage}[t]{0.55\textwidth}
            \centering
            \includegraphics[width=\textwidth]{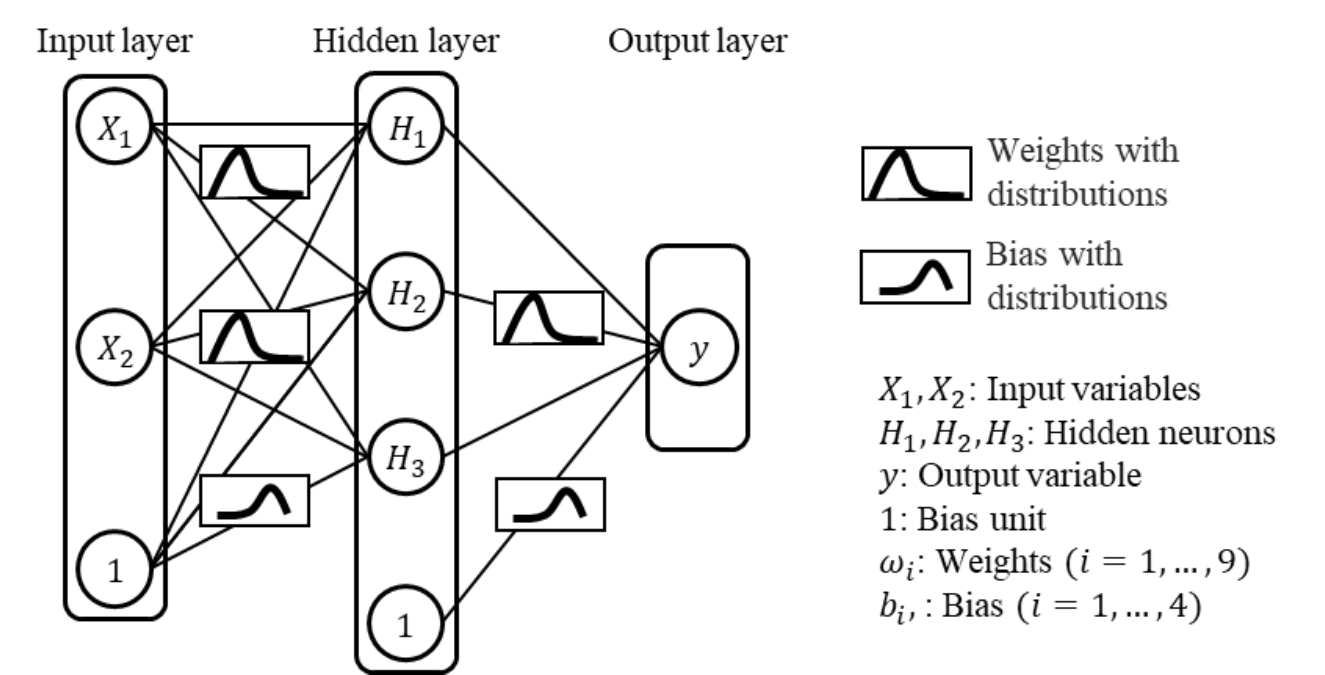}
            \caption*{(a)}
        \end{minipage}
        \hspace{0.0cm}
        \begin{minipage}[t]{0.34\textwidth}
            \centering
            \includegraphics[width=\textwidth]{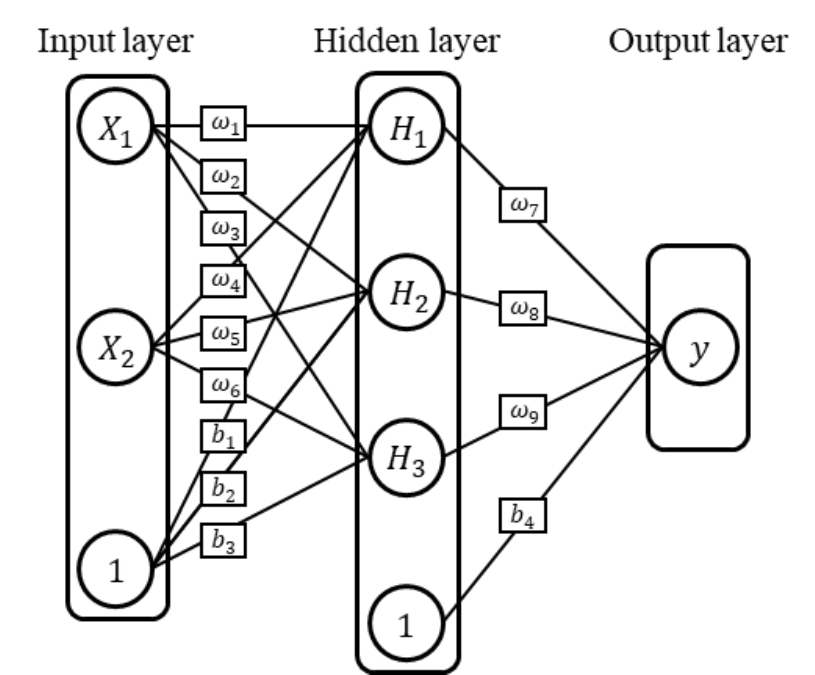}
            \caption*{(b)}
        \end{minipage}
        \caption{Schematic diagram of the neural network architecture. The network consists of input nodes $(X_1, X_2)$ and a bias unit denoted by $1$, which provides a constant input value. Also, $\boldsymbol{\omega}$ denotes the weights, and $\boldsymbol{b}$ represents the biases. Each hidden neuron $(H_1, H_2, H_3)$ applies a nonlinear activation function to the weighted sum of its inputs: (a) Bayesian neural networks; (b) Artificial neural network}
        \label{fig:BNN_structure}
    \end{figure}
    In contrast, BNN quantifies epistemic uncertainty by assigning probability distributions to the model parameters rather than treating them as deterministic values~\cite{goan2020bayesian}. 
    
    BNN, similar to GP, relies on Bayesian inference by treating the model parameters as a random vector $\boldsymbol{\Theta}$, which maps from $(\Omega, \mathcal{F})$ to a measurable space $(\mathbb{A}^{L_{\omega}+L_b}, \mathcal{B}^{L_{\omega}+L_b})$, with realizations $\boldsymbol{\theta}\in\mathbb{R}^{L_{\omega}+L_b}$. GP requires computational cost proportional to the cube of the data size due to covariance matrix inversion, making it impractical for large datasets. BNN, however, scales more efficiently through mini-batch stochastic optimization, where the training cost per iteration grows only linearly with data size. Moreover, while GP suffers from degraded kernel function performance in high-dimensional inputs, BNN effectively captures complex variable interactions through the hierarchical structure of neural networks.

    Following standard convention in Bayesian modeling, we use $p(\cdot)$ to denote prior, posterior, and likelihood distributions throughout this section. Given a training dataset $\mathcal{S}=\{\mathbf{x}^{(l)},y(\mathbf{x}^{(l)})\}_{l=1}^L$, we define the posterior distribution of the model parameters via Bayes' theorem, i.e.,
    \begin{align*}
        p(\boldsymbol{\theta} \vert \mathcal{S}) = \frac{p(\mathcal{S} \vert \boldsymbol{\theta})p(\boldsymbol{\theta})}{p(\mathcal{S})},
    \end{align*}
    where $p(\boldsymbol{\theta})$ is the prior distribution over the parameters and $p(\mathcal{S}|\boldsymbol{\theta})$ is the likelihood. The prior constraints the model parameters during training, effectively preventing overfitting without requiring a separate regularization term~\cite{mackay1992practical}. The denominator, $p(\mathcal{S})=\int p(\mathcal{S}|\boldsymbol{\theta})p(\boldsymbol{\theta})\mathrm{d}\boldsymbol{\theta}$, is called the evidence or marginal likelihood, function as a normalization constant indicating how well the model explains the data. However, in the context of deep neural networks, the high dimensionality of the weights makes the analytical calculation of the evidence through integration intractable. To address this challenge, variational inference serves as a standard approach~\cite{jordan1999introduction, graves2011practical}.
    
    Variational inference approximates the true posterior $p(\boldsymbol{\theta}|\mathcal{S})$ with a tractable distribution $q(\boldsymbol{\theta}|\boldsymbol{\phi})$ parameterized by the variational parameters $\boldsymbol{\phi}$. The objective is to minimize the Kullback-Leibler (KL) divergence between the approximate and true posterior distributions, i.e., $D_{KL}(q(\boldsymbol{\theta}|\boldsymbol{\phi})||p(\boldsymbol{\theta}))=\sum_{\theta} p(\boldsymbol{\theta}) \log \frac{p(\boldsymbol{\theta})}{q(\boldsymbol{\theta}|\boldsymbol{\phi})}$. Since directly minimizing this KL divergence requires the intractable true posterior, we equivalently maximize the evidence lower bound (ELBO) as
    \begin{align*}        \mathcal{L}_{\mathrm{ELBO}}=\mathbb{E}_{q(\boldsymbol{\theta}|\boldsymbol{\phi})}[\log p(\mathcal{S}|\boldsymbol{\theta})]-D_{KL}(q(\boldsymbol{\theta}|\boldsymbol{\phi})||p(\boldsymbol{\theta})).
    \end{align*}
    The first term measures the model's goodness of fit to the given data, and the second KL divergence term regularizes the approximate posterior to remain close to the prior. Maximizing the ELBO updates the variational parameters $\theta$, driving the approximate posterior $q(\boldsymbol{\theta}|\boldsymbol{\phi})$  toward the true posterior.
    
    We obtain the predictive distribution for the output $\tilde{y}_{\mathcal{N}}^{*}=\tilde{y}_{\mathcal{N}}(\mathbf{x}^*;\boldsymbol{\theta})$ at a new test input $\mathbf{x}^{*}$ by marginalizing over the posterior as 
    \begin{align*}
        p(\tilde{y}_{\cal{N}}^* | \mathbf{x}^{*}, \mathcal{S}) = \int p(\tilde{y}_{\cal{N}}^* | \mathbf{x}^{*}, \boldsymbol{\theta}) \cdot p(\boldsymbol{\theta}|\mathcal{S}) \mathrm{d} \boldsymbol{\theta}.
    \end{align*}
    The nonlinear activation functions and multi-layer structure of neural networks prevent a closed-form solution for this integral. We thus use Monte Carlo approximation by drawing $T$ weight samples $\boldsymbol{\theta}^{(t)}$ from the approximate posterior $q(\boldsymbol{\theta}|\boldsymbol{\phi})$  as 
    \begin{align*}
        p(\tilde{y}_{\cal{N}}^* | \mathbf{x}^{*}, \mathcal{S}) \approx \frac{1}{T} \sum_{t=1}^{T} p(\tilde{y}_{\cal{N}}^* | \mathbf{x}^*, \boldsymbol{\theta}^{(t)}).
    \end{align*}
    We then derive the predictive mean and variance as
    \begin{align}
        \mathbb{E}[\tilde{y}_{\cal{N}}(\mathbf{x}^{*}; \boldsymbol{\Theta})| \mathcal{S}] &= \tilde{\mu}_{\tilde{y}_{\cal{N}}}(\mathbf{x}^{*}; \boldsymbol{\theta}) \approx \frac{1}{T} \sum\limits_{t=1}^{T} \tilde{y}_{\cal{N}}(\mathbf{x}^{*}; \boldsymbol{\theta}^{(t)}),\\
        \mathbb{V}\mathrm{ar} [\tilde{y}_{\cal{N}}(\mathbf{x}^{*}; \boldsymbol{\Theta}) | \mathcal{S}] &=\tilde{\sigma}_{\tilde{y}_{\cal{N}}}^{2}(\mathbf{x}^{*}; \boldsymbol{\theta}) \approx \frac{1}{T} \sum\limits_{t=1}^{T} \tilde{y}_{\cal{N}}(\mathbf{x}^{*}; \boldsymbol{\theta}^{(t)})^{2} - \Big( \frac{1}{T} \sum\limits_{t=1}^{T} \tilde{y}_{\cal{N}}(\mathbf{x}^{*}; \boldsymbol{\theta}^{(t)}) \Big)^{2},
        \label{eq:bnn_var}
    \end{align}
    where $\tilde{y}_{\cal{N}}(\mathbf{x}^{*}; \boldsymbol{\theta}^{(t)})$ represents the neural networks output using the $t$-th weight sample. The predictive mean $\tilde{\mu}_{\tilde{y}_{\mathcal{N}}}$ represents the ensemble-average output, and the predictive variance $\tilde{\sigma}_{\tilde{y}_{\mathcal{N}}}$ quantifies the epistemic uncertainty, reflecting the model's confidence in its prediction. 
%

%%%%%%%%%%%%%%%%%%%%%%%%%%%%%%%%%%%%%%%%%%%
\section{Robust design optimization with BNN-PDD integration} \label{section:RDO for BNN-PDD}
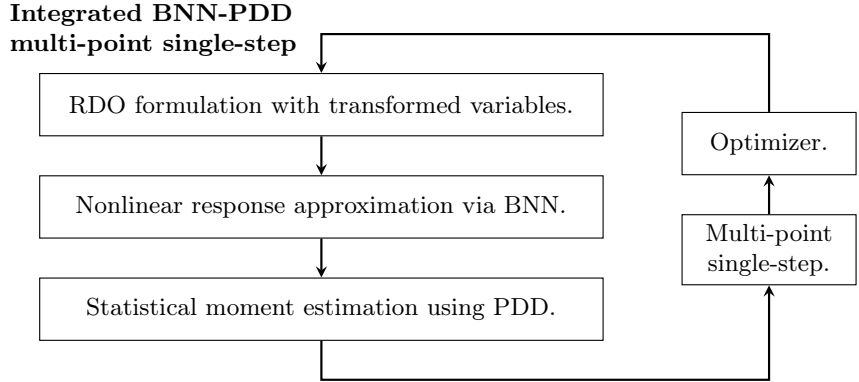
\begin{figure}[htb!]
    \centering
    \begin{minipage}[t]{0.20\textwidth}
        \vspace{0pt}
        \caption{The overall framework of the proposed BNN-PDD multi-point single-step optimization.}
        \label{fig:BNN_PDD_flowchart}
    \end{minipage} 
    \hfill
    \begin{minipage}[t]{0.75\textwidth}
        \vspace{0pt}
        \centering
        \resizebox{\linewidth}{!}{
            \begin{tikzpicture}[
                node distance=0.5cm,
                block/.style={rectangle, draw, fill=white, 
                    text width=7.0cm, align=center, minimum height=0.8cm, font=\small},
                arrow/.style={thick, ->, >=stealth},
                %container/.style={draw=green!40!black, thick, inner sep=0.4cm}
            ]
            
                \node (b1) [block] {RDO formulation with transformed variables.};
                \node (b2) [block, below=of b1] {Nonlinear response approximation via BNN.};
                \node (b3) [block, below=of b2] {Statistical moment estimation using PDD.};
            
                %\node (box) [container, fit=(b1) (b3)] {};
                
                \node [font=\bfseries\small, align=left, anchor=south west, xshift=-0.5cm, yshift=0.1cm] at (b1.north west) {Integrated BNN-PDD\\ multi-point single-step};
            
                % Optimizer
                \node (opt) [block, right=1.0cm of b1, yshift=-0.5cm, text width=2.0cm] {Optimizer.};
                \node (b4) [block, below=of opt, text width=2.0cm] {Multi-point single-step.};

                \draw [arrow] (b1) -- (b2);
                \draw [arrow] (b2) -- (b3);

                \draw [arrow] (b3.south)|- ($(b3.south) + (0, -0.5)$) -| (b4.south);
                \draw [arrow] (b4.north) -- (opt.south);
                \draw [arrow] (opt.north) |- ($(b1.north) + (0, 0.5)$) -| (b1.north);

                %\draw [arrow] (opt.north) |- ($(box.north) + (0, 0.5)$) -| (box.north);
            
                %\draw [arrow] (box.south) |- ($(box.south) + (0, -0.5)$) -| (opt.south);
    
            \end{tikzpicture}
        }
    \end{minipage}
\end{figure}

%%%%%%%%%%%%%%%%%%%%%%%%%%%%%%%%%%%%%%%%%%%
%
    This section proposes a novel method that combines efficient surrogate modeling with UQ to address the computational challenges of RDO. The proposed method leverages the predictive uncertainty from Bayesian inference to enable efficient exploration of the design space. Section~\ref{section:transformed random variables} formulates the RDO problem using transformed random variables. Section~\ref{section:Single step process} details the BNN-based single-step process, which integrates active learning for surrogate construction with the BNN-PDD method for statistical moment estimation. Section~\ref{section:MP} extends this method to a multi-point single-step process. Section~\ref{section:RDO_algorithm} outlines the complete algorithm of the proposed RDO framework. Figure~\ref{fig:BNN_PDD_flowchart} illustrates the overall procedure, comprising variable transformation, nonlinear response approximation via BNN, and moment estimation through PDD.

%

%%%%%%%%%%%%%%%%%%%%%%%%%%%%%%%%%%%%%%%%%%%
\subsection{Robust design optimization with transformed random variables} \label{section:transformed random variables}
%%%%%%%%%%%%%%%%%%%%%%%%%%%%%%%%%%%%%%%%%%%
%
The RDO estimates the mean and variance of system performance under uncertainty and derives an optimal design based on these statistics. During the optimization process, the probability distribution of $\mathbf{X}$ changes whenever design variables are updated. Since PDD basis functions must satisfy orthogonality with respect to the input PDF, this variation requires recalculating the basis functions at every iteration, resulting in substantial computational cost. 

To overcome this computational burden, we introduce transformed variables $\mathbf{U} = (U_1, \ldots, U_N)^{\intercal}$ that follow a fixed probability distribution independent of the design variables. Without loss of generality, we adopt the uniform distribution as a representative example. We assume each component of $\mathbf{U}$ follows a uniform distribution on a fixed interval $[\mathbf{a},\mathbf{b}]\subset \mathbb{R}^{N}$, where $\mathbf{a}=(a_1, \ldots, a_{N})^{T}$ and $\mathbf{b}=(b_1, \ldots, b_{N})^{\intercal}$ denote the lower and upper bound vectors, respectively. This formulation allows us to compute the PDD basis functions only once for the fixed uniform distribution, eliminating recalculation throughout the design iterations and significantly enhancing computational efficiency.
    
We map the transformed variables $\mathbf{U}$ to the original random variables $\mathbf{X}$ via the probability integral transformation. The joint CDF of $\mathbf{U}$ on $[\mathbf{a},\mathbf{b}]$ is
    \begin{align*}
        F_{\mathbf{U}}(\mathbf{u}) = \prod_{i=1}^{N} \frac{u_i - a_i}{b_i - a_i}, \quad \mathbf{u} \in [\mathbf{a},\mathbf{b}].
    \end{align*}
    We denote the CDF of $\mathbf{X}$ as $F_{\mathbf{X}}(\mathbf{x};\mathbf{d})$ as the CDF depends on the design vector $\mathbf{d}$. Equating $F_{\mathbf{X}}(\mathbf{x};\mathbf{d})$, we obtain the inverse mapping, i.e., 
    \begin{align}
        \mathbf{x} = F_{\mathbf{X}}^{-1} \left(\frac{\mathbf{u} - \mathbf{a}}{\mathbf{b} - \mathbf{a}};\mathbf{d}\right), 
        \label{eq:transform_variables}
    \end{align}
    where $F_{\mathbf{X}}^{-1}(\cdot)$ denotes the inverse CDF of $\mathbf{X}$. We note that the optimization still proceeds in the design variable space; the transformation only changes the domain over which we compute the statistical moments---from the design dependent $\mathbf{X}$ to the fixed $\mathbf{U}$ distribution. We formulate the RDO problem using the transformed functions $y_l(\mathbf{X})=y_l(F_{\mathbf{X}}^{-1}(\Phi(\mathbf{U});\mathbf{d}))=h_l(\mathbf{U};\mathbf{d})$, $l=0,1,\ldots,K$, as 
    \begin{equation}
        \begin{aligned}
            \min_{\mathbf{d} \in \mathcal{D} \subseteq \mathbb{R}^{M}} \quad& c_{0}(\mathbf{d}) := w_1 \frac{\mathbb{E}[h_{0}(\mathbf{U};\mathbf{d})]}{\mu_{0}} + w_2 \frac{\sqrt{\mathbb{V}\mathrm{ar}[h_{0}(\mathbf{U};\mathbf{d})]}}{\sigma_{0}}, \\
            \text{subject to} \quad& c_{l}(\mathbf{d}) := \alpha_l \sqrt{\mathbb{V}\mathrm{ar}[h_{l}(\mathbf{U};\mathbf{d})]} - \mathbb{E}[h_{l}(\mathbf{U};\mathbf{d})] \leq 0, \quad l=1, \ldots, K,\\
            &\qquad d_{k,L} \leq d_{k} \leq d_{k,U}, \quad k=1,\ldots, M.
        \end{aligned}
    \label{eq:RDO formulation}
    \end{equation} 
    This formulation enables the reuse of PDD basis function across all design iterations, significantly reducing computational cost while maintaining accuracy in statistical moment estimation.
% 

%%%%%%%%%%%%%%%%%%%%%%%%%%%%%%%%%%%%%%%%%%%
%
\subsection{BNN based single step process} \label{section:Single step process}
%%%%%%%%%%%%%%%%%%%%%%%%%%%%%%%%%%%%%%%%%%%
%
    This section presents an integrated framework for efficient statistical moment calculation in RDO. We first construct a BNN surrogate to approximate the system response, and then combine it with the transformed variable formulation and PDD for moment estimation.
%

%%%%%%%%%%%%%%%%%%%%%%%%%%%%%%%%%%%%%%%%%%%
%
\subsubsection{BNN modeling via active learning}\label{section:BNN_active learning}
%%%%%%%%%%%%%%%%%%%%%%%%%%%%%%%%%%%%%%%%%%%
%
\begin{figure}[htb!]
    \centering
    \begin{tikzpicture}[
        node distance=0.5cm,
        block1/.style={rectangle, draw, fill=white, 
            text width=4.0cm, align=center, minimum height=0.8cm, font=\small},
        block2/.style={rectangle, draw, fill=white, 
            text width=4.7cm, align=center, minimum height=0.8cm, font=\small},
        block3/.style={rectangle, draw, fill=white, 
            text width=1.5cm, align=center, minimum height=0.8cm, font=\small},
        block4/.style={rectangle, draw, fill=white, 
            text width=6.0cm, align=center, minimum height=0.8cm, font=\small},
        block5/.style={rectangle, draw, fill=white, 
            text width=4.0cm, align=center, minimum height=0.8cm, font=\small},
        arrow/.style={thick, ->, >=stealth},
        container/.style={draw=green!40!black, fill=green!10, rounded corners, 
            inner sep=0.2cm}
    ]
    
        \node (b1) [block1] {Initial training data.};
        \node (b2) [block1, below=1.0cm of b1] {Current training data $\mathcal{S}$.};
        \node (b3) [block2, below=0.8cm of b2] {Train BNN \& Evaluate error.};
        \node (b4) [block3, below=1.2cm of b3] {End.};

        \node (b5) [block4, right=0.8cm of b3, yshift=0.6cm] {Estimate uncertainty $\tilde{\sigma}^{2}_{\tilde{y}{\mathcal{N}}}(\mathbf{x}^{*(l)};\boldsymbol{\theta}),~l=1,\ldots,\bar{\bar{L}}$};
        \node (b6) [block1, above=of b5, minimum height=1.0cm] {Select top $N_{a}$ candidates\\ with uncertainty.};
        \node (b7) [block5, below=of b5] {Generate candidate samples $\mathcal{S}_{c}$.};

        \node [font=\bfseries\small, align=left, anchor=south east, xshift=0.0cm, yshift=0.2cm] at (b6.north east) {Active learning};
        
        \draw [arrow] (b1) -- (b2);
        \draw [arrow] (b2) -- (b3);
        \draw [arrow] (b3) -- (b4);
        \draw [arrow] (b3.east) -- ++(0.35cm, 0) |- (b5.west);
        \draw [arrow] (b5) -- (b6);
        \draw [arrow] (b6.west) -- ++(-0.8cm, 0) |- (b2.east);
        \draw [arrow] (b7) -- (b5);

        \begin{scope}[on background layer]
            \node [container, fit=(b2) (b3) (b5) (b6) (b7)] {};
        \end{scope}

        %\draw [arrow] (b3.south)|- ($(b3.south) + (0, -0.5)$) -| (b4.south);
        %\draw [arrow] (b4.north) -- (opt.south);
        %\draw [arrow] (opt.west) |- ($(b2.east) + (0, 0.5)$) -| (b2.east);

    \end{tikzpicture}
    \caption{Flowchart of the active learning algorithm using BNN.}
    \label{fig:active learning_flowchart}
    \end{figure}
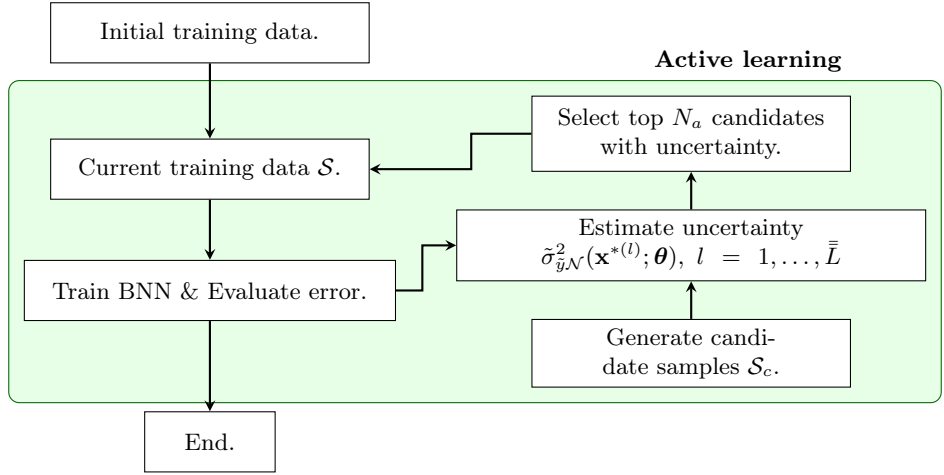
    We construct a BNN surrogate model within the physical design space to approximate the input-output relationship of the system. To acquire training data efficiently, we use Latin hypercube sampling (LHS)~\cite{helton2003latin, janssen2013monte}, which explores the design domain uniformly regardless of the probability distributions of design variables or transformed variables. We generate training samples via LHS and evaluate the corresponding response values from the original model to construct the initial dataset $\mathcal{S}$. The predictive performance of BNN (see Section~\ref{section:BNN}) depends on hyperparameter settings.  Instead of using fixed hyperparameters, we apply GP-based Bayesian optimization to determine the optimal hyperparameters for the dataset $\mathcal{S}$. With these optimal hyperparameters, BNN learns the posterior distribution of the weights and provides both predictive values and uncertainty estimates for arbitrary inputs.

    RDO requires maximizing surrogate model accuracy within a limited computational budget. Conventional random sampling, however, neglects model uncertainty, potentially leading to inefficient sampling placement. To address this, we use an uncertainty-based active learning strategy that leverages the epistemic uncertainty from BNN to reinforce regions of low model confidence. 
    
    Fig.~\ref{fig:active learning_flowchart} illustrates that the active learning algorithm operates through the following iterative procedure. We first generate $\bar{\bar{L}}$ candidate samples $\mathcal{S}_{c}=\{\mathbf{x}^{*(l)}\}_{l=1}^{\bar{\bar{L}}}$ through LHS within the currently defined subregion. We then evaluate the prediction variance $\tilde{\sigma}^{2}_{\tilde{y}{\mathcal{N}}}(\mathbf{x}^{*(l)};\boldsymbol{\theta})$ for each candidate $\mathbf{x}^{*(l)}$ using the trained BNN via \eqref{eq:bnn_var}. We sort all candidates in descending order of their predictive variance and select the top $N_a$ samples with the highest uncertainty, i.e., 
\begin{align*}
    \mathcal{S}_a = \left \{\mathbf{x}^{*(l)} \in \mathcal{S}_{c} \;\middle|\; \tilde{\sigma}^{2}_{\tilde{y}_{\mathcal{N}}}(\mathbf{x}^{*(l)};\boldsymbol{\theta}) \geq \tilde{\sigma}^{2}_{\tilde{y}_{\mathcal{N}}}(\mathbf{x}^{*(N_a)};\boldsymbol{\theta}),~l=1,\ldots,\bar{\bar{L}}>N_a \right\}.
 \end{align*}
    %\tilde{\sigma}^{2(l)}_{\tilde{y}{\mathcal{N}}}$ and select the top $N_{add}$ samples with the highest uncertainty. 
    Since these samples reside in regions where the model lacks confidence, we evaluate the original model at the selected samples to obtain actual response values and augment the training dataset accordingly. We then retain the BNN on the updated dataset, improving  prediction accuracy in previously uncertain regions. We access convergence by computing the relative error between the retrained model's predictions and the actual values for the newly added $N_a$ samples. We repeat this process until either the relative error falls below an acceptable threshold or the total number of training samples reaches a preset maximum. This active learning strategy enhances both computational efficiency and surrogate efficiency by concentrating function evaluations in regions of high uncertainty while avoiding unnecessary sampling elsewhere.  
    %We determine whether to terminate training based on prediction accuracy of the added samples. Specifically, we calculate the relative error between the retrained model's predictions and actual values, repeating the above process until the error of the newly added $N_{add}$ samples converges within the acceptable threshold or the total number of training samples reaches the preset maximum. Consequently, the proposed approach can enhance both the computational efficiency and predictive accuracy of the entire RDO process by minimizing function evaluations in unnecessary regions and refining the accuracy of the surrogate model.
%

%%%%%%%%%%%%%%%%%%%%%%%%%%%%%%%%%%%%%%%%%%%
%
\subsubsection{Statistical moment analysis via BNN-PDD} \label{section:BNN-PDD}
%%%%%%%%%%%%%%%%%%%%%%%%%%%%%%%%%%%%%%%%%%%
%
   \begin{figure}[htb!]
    \centering
    \begin{tikzpicture}[
        node distance=0.5cm,
        block1/.style={rectangle, draw, fill=white, 
            text width=5.7cm, align=center, minimum height=0.8cm, font=\small},
        block2/.style={rectangle, draw, fill=white, 
            text width=4.5cm, align=center, minimum height=0.8cm, font=\small},
        block3/.style={rectangle, draw, fill=white, 
            text width=4.5cm, align=center, minimum height=1.0cm, font=\small},
        block4/.style={rectangle, draw, fill=white, 
            text width=5.2cm, align=center, minimum height=1.0cm, font=\small},
        block5/.style={rectangle, draw, fill=white, 
            text width=5.2cm, align=center, minimum height=1.0cm, font=\small},
        arrow/.style={thick, ->, >=stealth}
    ]
    
        \node (b1) [block1] {Construct basis matrix $\mathbf{A}$ from $\mathbf{U}$.};
        \node (b2) [block2, below=of b1] {Transform $\mathbf{U}$ to $\mathbf{X}$ via~\eqref{eq:transform_variables}.};
        \node (b3) [block3, below=of b2] {Predict response vector $\mathbf{b}$\\ using trained BNN.};

        \node (b4) [block4, right=1.0cm of b2, yshift=0.0cm] {Compute expansion coefficient $\mathbf{c}$:\\ $\mathbf{c} = (\mathbf{A}^{\intercal} \mathbf{A})^{-1} \mathbf{A}^{\intercal} \mathbf{b}$.};
        \node (b5) [block5, right=1.0cm of b3] {Analytically derive the mean and variance using~\eqref{eq:pdd_mean_var}.};

        \draw [arrow] (b1) -- (b2);
        \draw [arrow] (b2) -- (b3);
        \draw [arrow] (b3.east) -- ++(0.35cm, 0) |- (b4.west);
        \draw [arrow] (b4) -- (b5);

    \end{tikzpicture}
    \caption{Flowchart of the statistical moment analysis using the proposed BNN-PDD method.}
    \label{fig:statistical moment_flowchart}
    \end{figure}
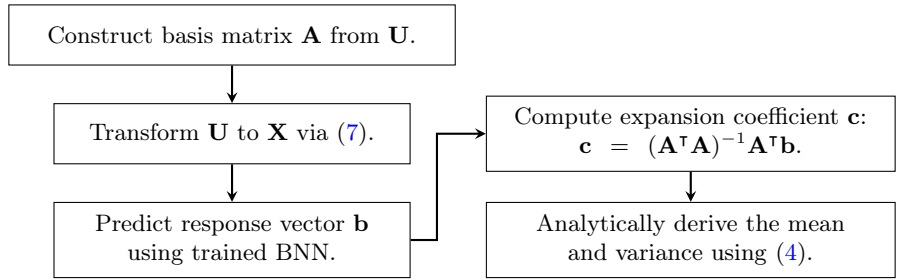
    
    We combines BNN for approximating nonlinear system responses with PDD for estimating statistical moments. When using BNN alone, statistical moment estimation relines on MCS, which suffers from slow convergence and numerical instability due to stochastic sampling noise. In contrast, PDD derives the mean and variance from expansion coefficients in closed-form expression, yielding noise-free, deterministic values throughout the optimization process. The proposed integration thus leverages the nonlinear approximation ability of BNN and the analytical efficiency of PDD, offering better numerical stability and computational performance compared to MCS.

    We apply the transformed variable formulation from Section~\ref{section:transformed random variables} to eliminate the cost of recomputing orthonormal polynomials whenever design variables change. To maximize the estimation accuracy of PDD expansion coefficients, we use the Sobol sequence---a low-discrepancy sequence that distributes samples uniformly even in high-dimensional spaces. This quasi-random sampling satisfies the orthogonality of PDD basis functions with fewer samples and accelerates convergence compared to simple random sampling.

    Figure~\ref{fig:statistical moment_flowchart} illustrates the overall computational procedure. We first precompute and fix the basis matrix $\mathbf{A}$ using Legendre polynomials evaluated at the Sobol-generated $\mathbf{U}$ samples; this calculation occurs only once before optimization begins. 
    %for the $\mathbf{U}$ generated by the Sobol sequence, performing this calculation only once before optimization. 
    At each optimization iteration, given a design vector $\mathbf{d}$, we convert $\mathbf{U}$ to $\mathbf{X}$ through the inverse transformation in~\eqref{eq:transform_variables} and feed the resulting $\mathbf{X}$ into the trained BNN to predict the response vector $\mathbf{b}$. We then compute the expansion coefficient vector as $\mathbf{c} = (\mathbf{A}^{\intercal} \mathbf{A})^{-1} \mathbf{A}^{\intercal} \mathbf{b}$ via least squares. Exploiting the orthonormality of the basis functions, we analytically obtain the mean and variance from $\mathbf{c}$ using~\eqref{eq:pdd_mean_var} without additional simulations, and directly evaluate the RDO objective and constraint functions in~\eqref{eq:weighted sum}. This procedure reduces statistical moment estimation to a sequence of variable transformation and matrix operations, completely avoiding basis function reconstruction across design iterations.
%   using the matrix $\mathbf{A}$ and the updated vector $\mathbf{b}$. The obtained $\mathbf{c}$ convert analytically into the mean and variance using~\eqref{eq:pdd_mean_var}, respectively, through the orthonormality of basis functions without additional simulations, directly utilizing them for the RDO objective function defined in~\eqref{eq:weighted sum}. This approach enables the estimation of statistical moments through variable transformation and matrix operations without complex reconstruction of basis functions.

%
%%%%%%%%%%%%%%%%%%%%%%%%%%%%%%%%%%%%%%%%%%%
\subsection{Multi-point single-step process via BDD-PDD} \label{section:MP}
%%%%%%%%%%%%%%%%%%%%%%%%%%%%%%%%%%%%%%%%%%%
%
The single-step method in Section~\ref{section:Single step process} constructs a BNN surrogate over the entire design space $\mathcal{D}$. For large design spaces or highly nonlinear responses, however, this global approach demands high-order basis functions or numerous samples for accurate approximation, degrading computational efficiency. A lower-accuracy approximation resulting from limited training samples may also fail to locate the true global optimum.

To address these challenges, we adopt a multi-point single-step (MPSS) strategy that decomposes the original RDO problem into a series of local RDO subproblems, where the objective and constraint functions in each subproblem are approximated within the corresponding local region~\cite{toropov1993multiparameter}.
    %The single-step method in Section~\ref{section:Single step process} represents a global approach that constructs a BNN surrogate for the entire design space. However, for large design spaces or highly nonlinear response functions, this global approach necessitates high-order basis functions or numerous samples for accurate approximation, potentially degrading computational efficiency. Using a lower accurate approximation induced by limited training samples may fail to find the accurate global optimum. To address these challenges, we use a multi-point single-step strategy. According to this method, the original RDO problem is replaced by a series of simplified local RDO sub-problems, where the objective and constraint functions of each sub-problem represent multi-point approximations~\cite{toropov1993multiparameter} in the corresponding local region.
    \begin{figure}[htb!]
        \centering
        \includegraphics[width=0.6\textwidth]
        {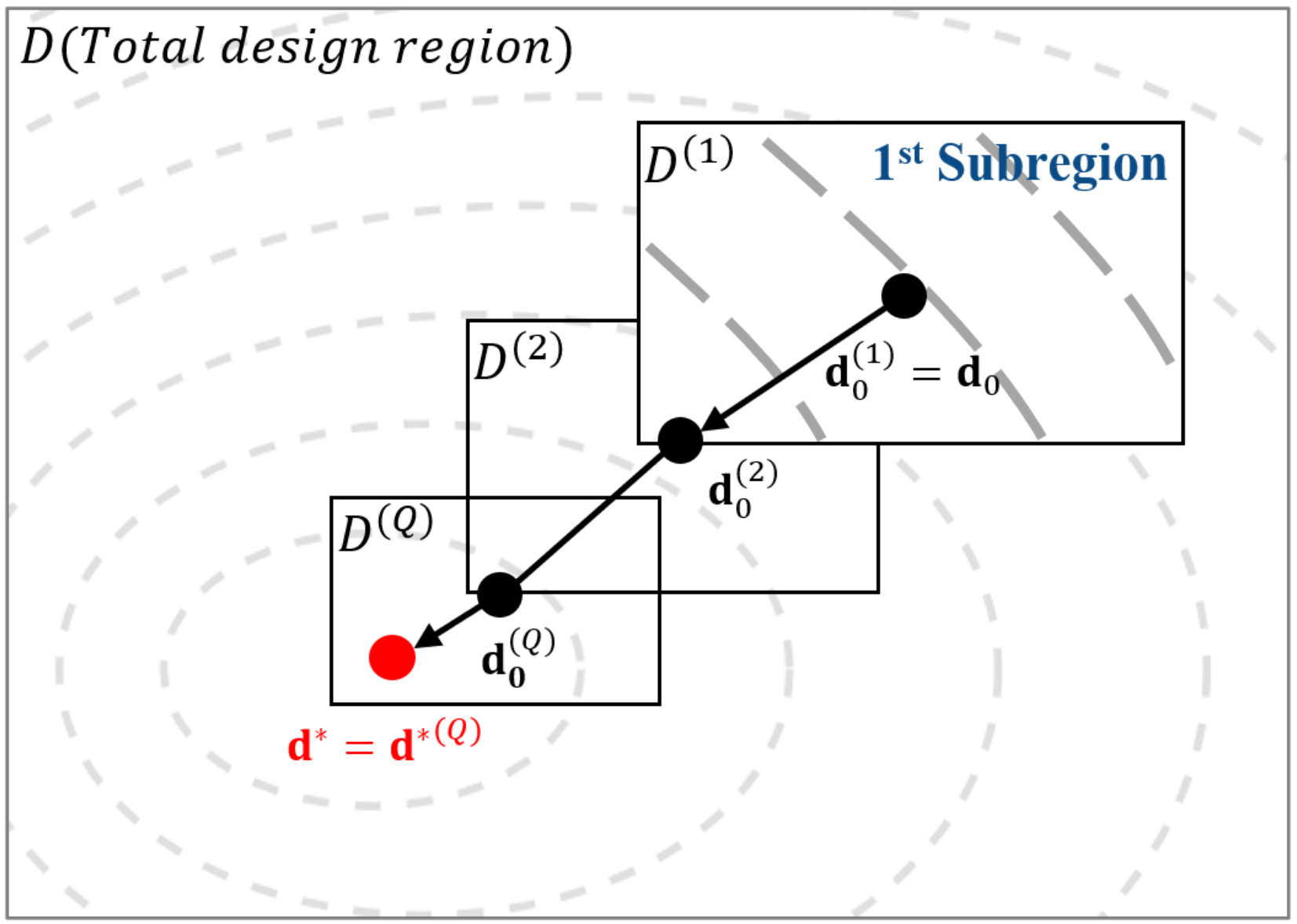}
    \caption{Schematic description of the multi-point single-step process.}
    \label{fig:MP_process}
    \end{figure}
    
    Figure~\ref{fig:MP_process} illustrates the MPSS process. At the $q$-th iteration, we define a subregion $\mathcal{D}^{(q)}$ centered on the current initial design point $\mathbf{d}_{0}^{(q)}$ within the global design space $\mathcal{D} = \prod_{k=1}^{M} [d_{k,L}, d_{k,U}]$ as 
    \begin{align*}
        \mathcal{D}^{(q)} = \prod_{k=1}^{M} \Bigg{[} d_{k,0}^{(q)} - \beta_{k}^{(q)} \frac{(d_{k,U}-d_{k,L})}{2}, \, d_{k,0}^{(q)} + \beta_{k}^{(q)} \frac{(d_{k,U}-d_{k,L})}{2} \Bigg{]}, \quad q=1, \ldots, Q.
    \end{align*}
    Here, $0 < \beta_{k}^{(q)} \leq 1$ is the sizing parameter controlling  the subregion extent for the $k$-th design variable. We solve each subproblem using low-order BNN-PDD approximations within the defined subregion. The optimal design solution derived from the current subproblem then serves as the initial design point for the next iteration. We update $\beta_k^{(q)}$ and repeat this process until convergence to the optimum. 
%

%%%%%%%%%%%%%%%%%%%%%%%%%%%%%%%%%%%%%%%%%%%
\subsection{The complete algorithm for robust design optimization} \label{section:RDO_algorithm}
%%%%%%%%%%%%%%%%%%%%%%%%%%%%%%%%%%%%%%%%%%%
%
\begin{figure}[htb!]
    \centering
    \begin{minipage}[t]{0.20\textwidth}
        \vspace{0pt}
        \caption{Flowchart of the multi-point single-step BNN-PDD RDO algorithm.}
        \label{MP_flowchart}
    \end{minipage}
    \hfill
    \begin{minipage}[t]{0.75\textwidth}
        \vspace{0pt}
        \centering
        \resizebox{\linewidth}{!}{
            \begin{tikzpicture}[
                node distance=0.6cm and 1cm,
                process/.style={rectangle, draw=black, thick, minimum width=8cm, minimum height=1cm, align=center, fill=white},
                decision/.style={diamond, draw=black, thick, aspect=2, minimum width=3cm, align=center, fill=white},
                io/.style={rectangle, draw=black, thick, minimum width=3cm, minimum height=1cm, align=center, fill=white},
                arrow/.style={thick, ->, >=stealth},
                subproc/.style={rectangle, draw=black, thick, minimum width=5cm, minimum height=1cm, align=center, fill=white},
                sideproc/.style={rectangle, draw=black, thick, minimum width=3cm, minimum height=1cm, align=center, fill=white}]

                % Step 1 ~ 3
                \node (s1) [process] {Step 1: Identify initial design; $\mathbf{d}_0$.};
                \node (s2) [process, below=of s1] {Step 2: Initialization; $q=1$, set $\mathbf{d}_0^{(q)} = \mathbf{d}_0$.};
                \node (s3) [process, below=of s2] {Step 3: Transform $\mathbf{X}$ to $\mathbf{U}$ in Section.~\ref{section:transformed random variables}.};
                
                % Step 4
                \node (s4_1) [subproc, below=1.5cm of s3, xshift=-2.25cm] {Step 4-1: \\Calculate $\rho_{k}^{(q)} := \frac{|d_{k,0}^{(q)} - d_{k,0}^{(q-1)}|}{(d_{k,U} - d_{k,L})/2}$.};
                \node(s4_2) [subproc, below=0.6cm of s4_1] {Step 4-2: \\If $\rho_{k}^{(q)} \leq \epsilon_1$.};
                \node (s4_dec) [sideproc, right=1.5cm of s4_2] {Decrease $\beta_k^{(q)}$.};
                % Step 4-3
                \node (s4_3) [subproc, below=0.6cm of s4_2] {Step 4-3: \\ If $\epsilon_1 < \rho_k^{(q)} \leq \epsilon_2$.};
                \node (s4_main) [sideproc, right=1.5cm of s4_3] {Maintain $\beta_k^{(q)}$.};
                % Step 4-Else
                \node (s4_inc) [subproc, below=0.6cm of s4_3] {Increase $\beta_k^{(q)}$.};
                % Step 4 End
                \node (s4_end) [rectangle, draw=black, thick, minimum width=3cm, align=center, below=of s4_inc] {Go to Step 5.};
                % Draw Orange Box around Step 4
                \node[draw=orange!100, thick, fit=(s4_1) (s4_dec) (s4_main) (s4_inc) (s4_end), inner sep=20pt] (box4) {};
                
                % Step 5 ~ 8
                \node (s5) [process, below=1.0cm of box4] {Step 5: At $D^{(q)}$, train BNN surrogate using \\ sample reuse and active learning.};
                \node (s6) [process, below=of s5] {Step 6: Estimate statistical moments using PDD.};
                \node (s7) [process, below=of s6] {Step 7: Solve $q$-th local RDO problem.};
                \node (s8) [decision, below=of s7] {Step 8: \\ Converge?};
                \node (end) [io, below=of s8] {End; $\mathbf{d}^*$.};
                \node (update) [rectangle, draw=black, thick, minimum width=2cm, minimum height=1cm, right=2.5cm of s5] {$q = q + 1$.};
                
                % --- Connections (Arrows) ---
                % Step 1-3
                \draw [arrow] (s1) -- (s2);
                \draw [arrow] (s2) -- (s3);
                \draw [arrow] (s3) -- (box4.north); % Step 3 -> Box Top
                % Step 4
                \draw [arrow] (s4_1) -- (s4_2);
                \draw [arrow] (s4_2) -- node[anchor=south] {Yes} (s4_dec);
                \draw [arrow] (s4_2) -- node[anchor=west] {No} (s4_3);
                \draw [arrow] (s4_3) -- node[anchor=south] {Yes} (s4_main);
                \draw [arrow] (s4_3) -- node[anchor=west] {No} (s4_inc);
                \draw [arrow] (s4_inc) -- (s4_end);
                % Side logic connections in Step 4
                \draw [arrow] (s4_dec.east) -- ++(0.3,0) |- (s4_end.east);
                \draw [arrow] (s4_main.east) -- ++(0.3,0) |- (s4_end.east);
                % Step 5-8
                \draw [arrow] (box4.south) -- (s5); % Box Bottom -> Step 5
                \draw [arrow] (s5) -- (s6);
                \draw [arrow] (s6) -- (s7);
                \draw [arrow] (s7) -- (s8);
                \draw [arrow] (s8) -- node[anchor=east] {Yes} (end);
                % Loop Back
                \draw [arrow] (s8.east) -| node[anchor=south, xshift=-5.5cm] {No} (update.south);
                \draw [arrow] (update.north) |- (s3.east);
            \end{tikzpicture}
        }
    \end{minipage}
\end{figure}

    Figure~\ref{MP_flowchart} illustrates the complete procedure of the proposed RDO algorithm. We detail each step below.
\begin{description}[leftmargin=1.6cm, style=sameline, font=\normalfont]
    \item[Step 1:]
    We determine an initial design vector $\mathbf{d}_{0}$ by first generating training data through LHS over the entire design space $\mathcal{D}$ and constructing an initial BNN model that captures the global trend. We apply the active learning strategy from Section~\ref{section:BNN_active learning} to ensure the reliability of this initial model. We then perform an approximate global optimization using~\eqref{eq:RDO formulation}  and set the resulting solution as the starting point $\mathbf{d}_0$. Go to Step~2.
    \item[Step 2:]
    We set the subregion iteration index $q = 1$ and designate $\mathbf{d}_{0}^{(q)} = \mathbf{d}_{0}$ as the center of the first subregion. We also initialize the subregion parameters required for algorithm execution, including the subregion sizing parameter $\beta_{k}^{(q)}$ and convergence tolerances. Go to Step~3.
    \item[Step 3:]
    We apply the variable transformation in ~\eqref{eq:transform_variables}, mapping the input random variable $\mathbf{X}$ to the transformed variable $\mathbf{U}$ that follows a fixed uniform distribution. Go to Step~4.
    \item[Step 4:]
    We compute the normalized displacement
\begin{align*}
    \rho_{k}^{(q)} := \frac{\left| d_{k,0}^{(q)} - d_{k,0}^{(q-1)} \right|}{(d_{k,U} - d_{k,L})/2},
\end{align*}
defined as the ratio of the center point movement to half the global design range. We then update the sizing parameter $\beta_{k}^{(q)}$ according to the following rule:
\begin{itemize}
    \item[] if $\rho_{k}^{(q)} \leq \epsilon_{1}$: decrease $\beta_{k}^{(q)}$ (solution approaching convergence; contract the subregion);
    \item[] if $\epsilon_{1} < \rho_{k}^{(q)} \leq \epsilon_{2}$: maintain $\beta_{k}^{(q)}$ (search range adequate);
    \item[] otherwise: increase $\beta_{k}^{(q)}$ (solution still exploring; expand the subregion).
\end{itemize}
This adaptive sizing enables the algorithm to begin with a broad search region and progressively narrow it as the solution converges toward the optimum. Go to Step~5.
    %We determine the sizing parameter $\beta_{k}^{(q)}$ at the current $q$-th iteration based on the rate of change between the current initial design $\mathbf{d}_{0}^{(q)}$ and the optimal design $\mathbf{d}^{*(q)}$. The rate of change $\rho_{k}^{(q)} := |d_{k,0}^{(q)} - d_{k,0}^{(q-1)}| / [(d_{k,U} - d_{k,L})/2]$ is defined as the ratio of the center point's movement distance relative to the local design space. If $\rho_{k}^{(q)} \leq \epsilon_{1}$, we judge that the solution has entered the convergence region and decrease $\beta_{k}^{(q)}$ to contract the subregion. If $\epsilon_{1} < \rho_{k}^{(q)} \leq \epsilon_{2}$, we consider the search range appropriate and the maintain $\beta_{k}^{(q)}$. Otherwise, increase $\beta_{k}^{(q)}$ and go to Step 5. This enables the optimization to start with a broad region in the early stages and progressively narrow the area for refined search as it converges toward the optimum.
    \item[Step 5:]
    We train the BNN surrogate within the current subregion $D^{(q)}$ using two strategies for data efficiency. First, we identify the overlap between the current subregion $\mathcal{D}^{(q)}$ and the previous one $\mathcal{D}^{(q-1)}$, and directly reuse existing training samples within this region without additional model evaluations. Second, we apply the uncertainty-based active learning from Section~\ref{section:BNN_active learning} to add new samples in regions of high predictive uncertainty, improving surrogate prediction accuracy where it is most needed. Go to Step~6.
    \item[Step 6:]
    We analytically compute the mean and variance of the response function via~\eqref{eq:pdd_mean_var} using the trained BNN model and the precomputed PDD basis functions.  Go to Step~7.
    \item[Step 7:]
    We solve the local RDO problem within $\mathcal{D}^{(q)}$ using the estimated statistical moments, obtaining the optimal solution $\mathbf{d}^{*(q)}$. Go to Step~8.
    \item[Step 8:]
    We evaluate convergence by checking the following criteria between the current optimal design $\mathbf{d}^{*(q)}$ and the current initial design $\mathbf{d}_{0}^{(q)}$:  if $||\mathbf{d}^{*(q)} - \mathbf{d}_{0}^{(q)}|| \leq \epsilon_3$ or $||c_{0}(\mathbf{d}^{*(q)}) - c_{0}(\mathbf{d}_{0}^{(q)})|| \leq \epsilon_4$, the algorithm terminates and outputs $\mathbf{d}^{*}=\mathbf{d}^{*(q)}$ as the final optimum. Otherwise, we set $\mathbf{d}_{0}^{(q+1)} = \mathbf{d}^{*(q)}$, increment $q = q + 1$, and return to Step 3.
\end{description}
%

%%%%%%%%%%%%%%%%%%%%%%%%%%%%%%%%%%%%%%%%%%%
\section{Numerical examples} \label{section:Numerical examples}
%%%%%%%%%%%%%%%%%%%%%%%%%%%%%%%%%%%%%%%%%%%
%
    This section presents two numerical examples to validate the effectiveness and computational efficiency of the proposed BNN-PDD multi-point single-step RDO framework. Section~\ref{section:Rastrigin function} uses a benchmark mathematical RDO problem to evaluate the global search ability and convergence stability of the proposed algorithm. Section~\ref{section:SPMSM} demonstrates applicability to high-dimensional problems through a shape RDO for an electric motor.
    
    In Example~1, we compare the proposed method against Gaussian process (GP) regression using the same number of training samples for objective performance evaluation. GP is a nonparametric Bayesian model that assumes function values follow a multivariate normal distribution, defining a zero-mean prior $y(\mathbf{x}) \sim \mathcal{GP}(0, k(\mathbf{x}, \mathbf{x}'))$ for the latent function $y(\mathbf{x})$, where $k(\mathbf{x}, \mathbf{x}')$ denotes the covariance function. We use the Mat\'{e}rn kernel $k(\mathbf{x}, \mathbf{x}') = \sigma^2 \frac{2^{1-\nu}}{\Gamma (\nu)} (\sqrt{2 \nu} \frac{r}{\ell})^{\nu}K_{\nu}(\sqrt{2 \nu} \frac{r}{\ell})$ to control function smoothness. Here, $\sigma^2>0$ is the signal variance controlling the overall amplitude of the function,  $r= \| \mathbf{x}-\mathbf{x}' \|$ is the Euclidean distance between two points, $\Gamma(\cdot)$ is the gamma function, $K_{\nu}(\cdot)$ is the modified Bessel function of the second kind, $\ell>0$ is the length scale controlling the rate of variation, and $\nu>0$ is the smoothness parameter determining the differentiability order of the function.
    
    In both examples, we formulate the problem as multi-objective optimization that simultaneously minimizes the mean and variance of the responses, and apply the weighted sum method to convert the multi-objective into a single objective function. For all cases, we set the PDD truncation parameters as polynomial order $m=3$ and interaction order $S=1$. Convergence tests confirmed that moment estimates exhibit negligible changes beyond these orders. We define the tolerance and subregion control parameters for MPSS algorithm as $\epsilon_1=0.1, \epsilon_2=0.5$, and $\epsilon_3=\epsilon_4=0.01$, with an initial subregion sizing parameter of $\beta_{k}^{(1)}=0.15$. For the local optimizer at each step, we adopt the differential evolution (DE) algorithm, which provides effective gradient-free global search for multi-modal objective functions. Comparative test with genetic algorithm (GA) showed that DE achieves faster convergence and better solutions quality. We obtain all numerical results using Python 3.11 on an Intel i5-13400 processor with 32 GB of RAM. The finite element analysis (FEA) in Example~2 was performed using JMAG (version 23.2).
%

%%%%%%%%%%%%%%%%%%%%%%%%%%%%%%%%%%%%%%%%%%%
%
\subsection{Example 1: RDO of a benchmark mathematical function} \label{section:Rastrigin function}
    We validate the performance of the proposed RDO algorithm through a mathematical benchmark. We compare three approaches: (1) the proposed BNN-PDD MPSS method, (2) GP-PDD MPSS which replaces BNN with GP as the surrogate model under the same domain settings and active learning strategy, and (3) MCS, which directly performs MCS over the global design space without a surrogate model. We consider two cases: a low-dimensional setting (Case 1: $N=2$) and a high-dimensional setting (Case 2: $N=10$). 
%

%%%%%%%%%%%%%%%%%%%%%%%%%%%%%%%%%%%%%%%%%%%
%
\subsubsection{Problem definition} \label{example1_problem definition}
    We model each random variable $X_{i}$ as a truncated Gaussian distribution with mean $d_{i}$ and standard deviation $\sigma_{i}$. We determine $\sigma_{i}$ applying the $3 \sigma$ rule to a manufacturing tolerance of $\pm 0.01$, yielding $\sigma_{i} = 0.01/3$. We formulate the RDO problem to minimize the weighted sum of the mean and standard deviation as defined in~\eqref{eq:weighted sum}:
    \begin{equation*}
        \begin{aligned}
            \min\limits_{\mathbf{d} \in \mathcal{D}} \quad c_{0}(\mathbf{d}) =&  w_{1} \frac{\mathbb{E}_{\mathbf{d}}[y(\mathbf{X})]}{\mathbb{E}_{\mathbf{d}_{0}}[y(\mathbf{X})]} + w_{2} \frac{\sqrt{\mathbb{V}\mathrm{ar}_{\mathbf{d}}[y(\mathbf{X})]}}{\sqrt{\mathbb{V}\mathrm{ar}_{\mathbf{d}_{0}}[y(\mathbf{X})]}}, \\
            -5.12 \leq& d_{i} \leq 5.12, \quad i=1,\ldots, M,
        \end{aligned}
    \end{equation*}
    where $y(\mathbf{X})$ is the Rastrigin function:
    \begin{align*}
        y(\mathbf{X}) = 10N + \sum\limits_{i=1}^{N}[X_{i}^{2} - 10 \cos(2 \pi X_{i})].
    \end{align*}
    We set equal weight coefficients $w_{1}=w_{2}=0.5$. The terms $\mathbb{E}_{\mathbf{d}_{0}}[y(\mathbf{X})]$ and $\sqrt{\mathbb{V}\mathrm{ar}_{\mathbf{d}_{0}}[y(\mathbf{X})]}$ denote the mean and standard deviation at the initial design $\mathbf{d}_{0}$, respectively, serving as normalization factors. 
    
\begin{figure}[htb!]
        \centering
        \includegraphics[width=0.8\textwidth, trim=0 280 0 270]{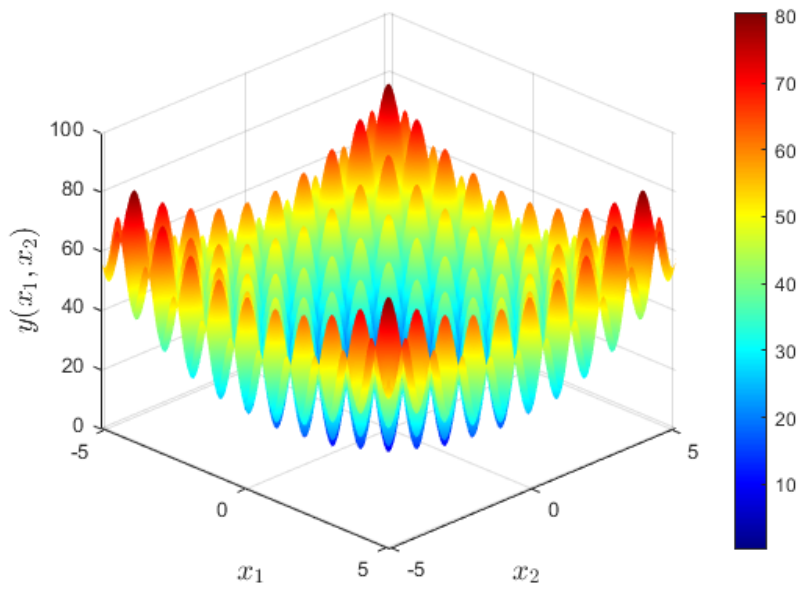}
    \caption{Visualization of the two-dimensional Rastrigin function. The highly multimodal landscape with numerous local minima highlights the difficulty of locating the global optimum at the origin.}
    \label{fig:Rastrigin_surface}
    \end{figure}

    The Rastrigin function is a widely used multi-modal benchmark with numerous local minima, posing a significant challenge for locating the global optimum. Fig.~\ref{fig:Rastrigin_surface} presents that this function has a unique global minimum at the origin $(0,0)^{\intercal}$.
\subsubsection{Method} \label{example1_method}
 We perform RDO for both the two-dimensional (2D) problem (Case 1) and the ten-dimensional (10D) problem (Case 2) problems using three methods: BNN-PDD MPSS, GP-PDD MPSS, and direct MCS. For fair comparison, we train both the BNN and GP surrogates using 500 training samples identically generated through LHS in each subregion. 

%%%%%%%%%%%%%%%%%%%%%%%%%%%%%%%%%%%%%%%%%%%
%
\subsubsection{Results for Case 1 ($N=M=2$)} \label{example 1_result_case1}
\begin{figure}[htb!]
    \centering
    \begin{minipage}[t]{0.48\textwidth}
        \centering
        \includegraphics[width=\textwidth]{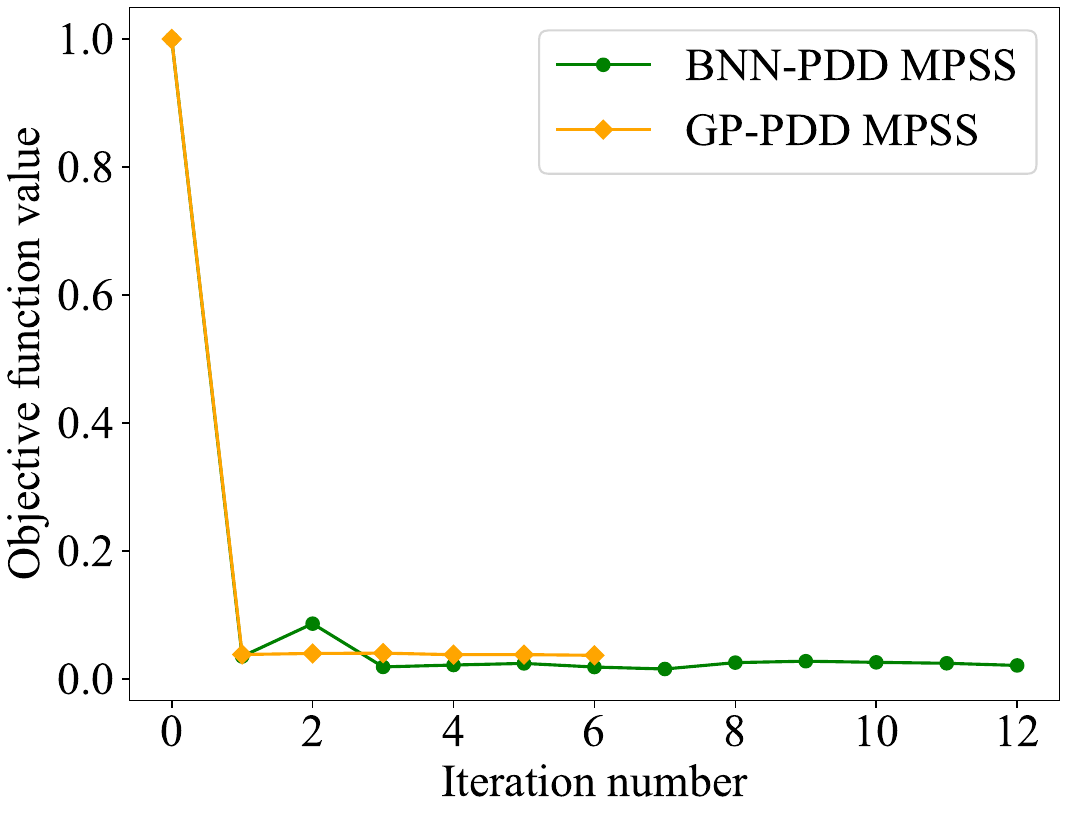}
        \caption*{(a) Case 1 (2D) - BNN-PDD MPSS \& GP-PDD MPSS}
    \end{minipage}
    \hfill
    \begin{minipage}[t]{0.48\textwidth}
        \centering
        \includegraphics[width=\textwidth]{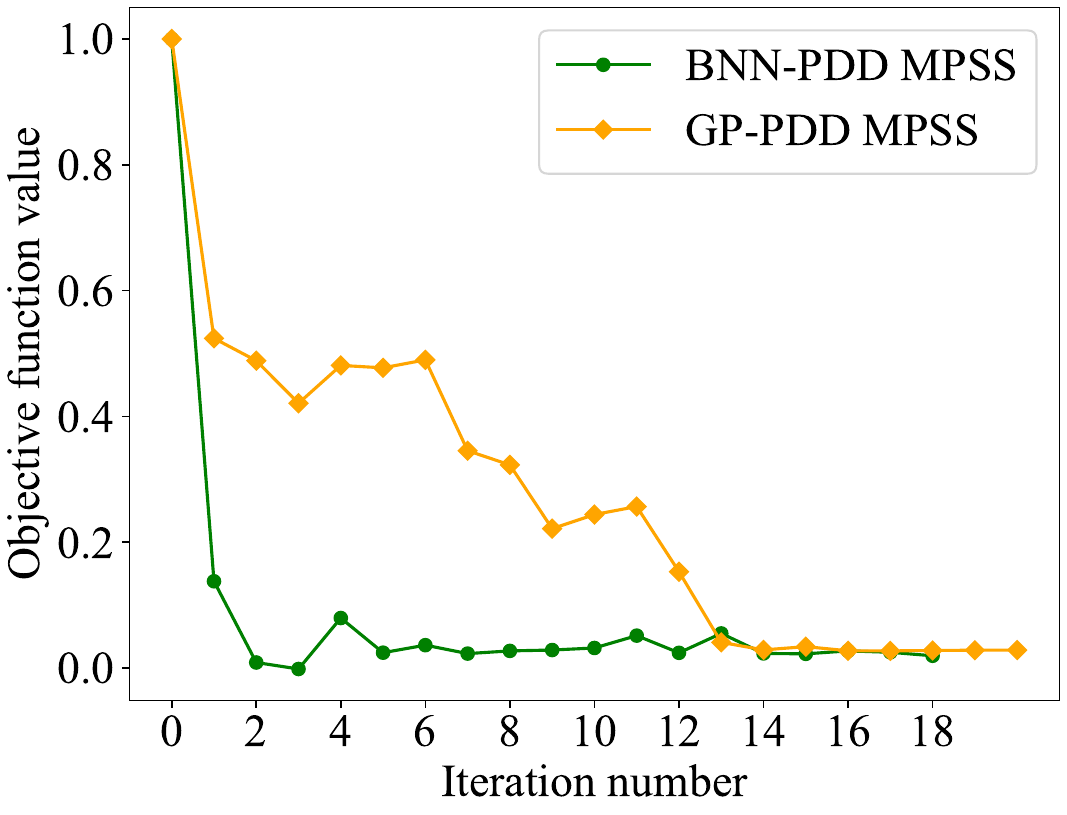}
        \caption*{(b) Case 2 (10D) - BNN-PDD MPSS \& GP-PDD MPSS}
    \end{minipage}
    \caption{Optimization convergence history of the objective function: (a) Case 1 using BNN-PDD MPSS and GP-PDD MPSS, (b) Case 2 using BNN-PDD MPSS and GP-PDD MPSS.}
    \label{fig:Rastrigin_obj}
    \end{figure}

    Figure~\ref{fig:Rastrigin_obj} shows the convergence history of the objective function for both Cases 1 and 2. The horizontal axis represents the subregion iteration index, and the vertical axis indicates the corresponding objective function value. Since no initial design point is prescribed, we use $\mathbf{d}^{(0)}$ obtained in Step 1 of Section.~\ref{section:RDO_algorithm} as the initial design; consequently, the initial values for BNN and GP differ. As shown in Fig.~\ref{fig:Rastrigin_obj}(a), both BNN-PDD MPSS and GP-PDD MPSS converge stably, with the objective function decreasing monotonically over iterations. In Case 1, both the BNN and GP-based approaches reaches the vicinity of the global optimum.
    
\begin{figure}[htb!]
    \centering
    \begin{minipage}[t]{0.48\textwidth}
        \centering
        \includegraphics[width=\textwidth]{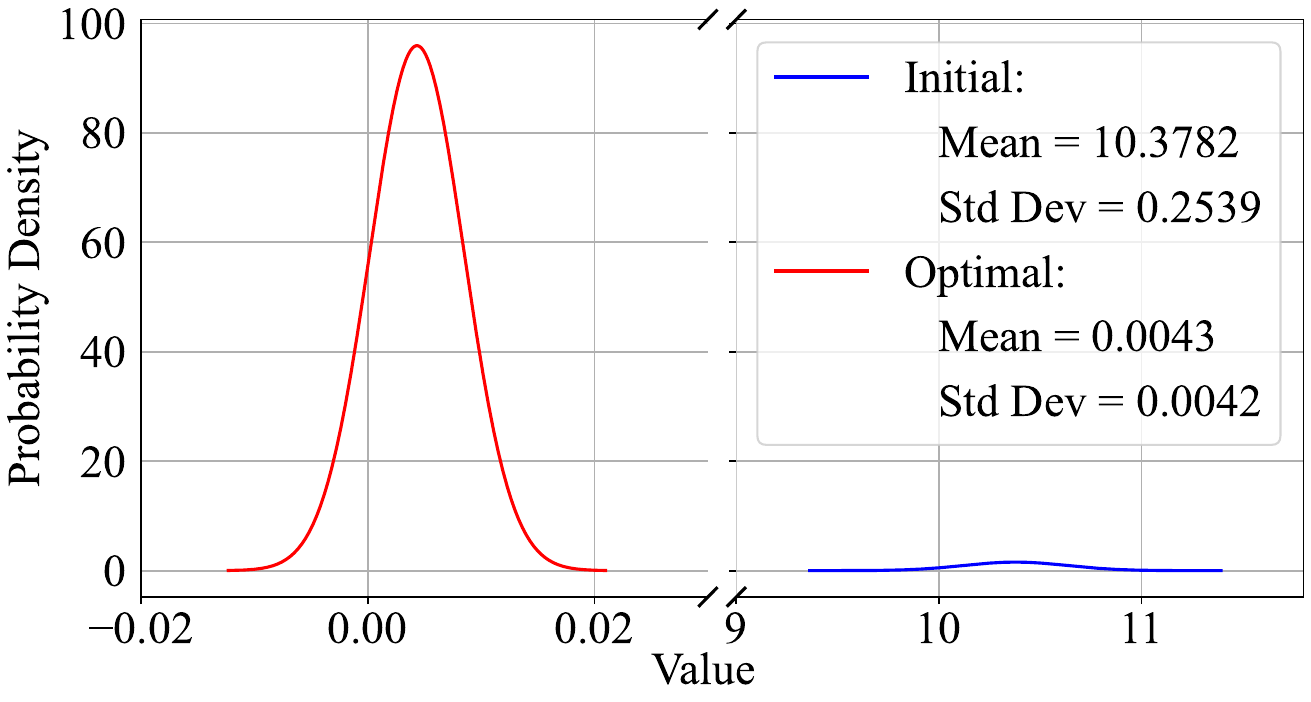}
        \caption*{(a) Case 1 (2D) - BNN-PDD MPSS}
    \end{minipage}
    \hfill
    \begin{minipage}[t]{0.48\textwidth}
        \centering
        \includegraphics[width=\textwidth]{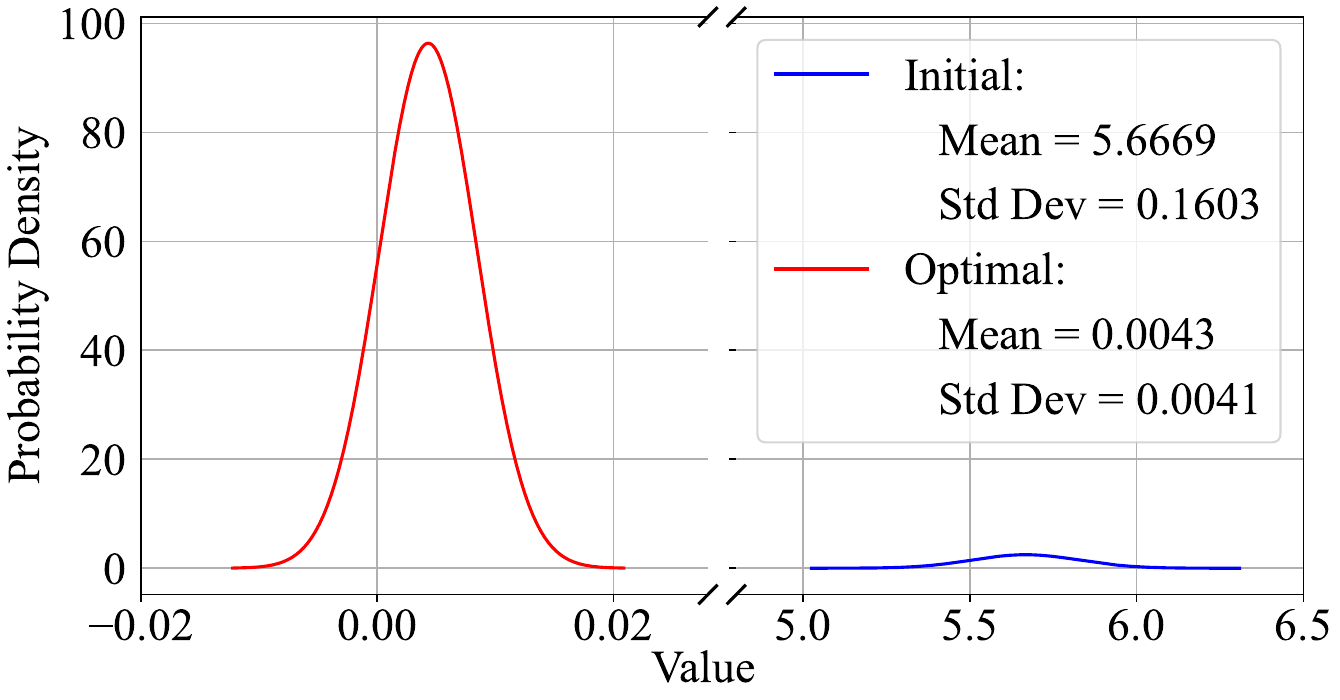}
        \caption*{(b) Case 1 (2D) - GP-PDD MPSS}
    \end{minipage}
    \\[0.5cm]
    \begin{minipage}[t]{0.48\textwidth}
        \centering
        \includegraphics[width=\textwidth]{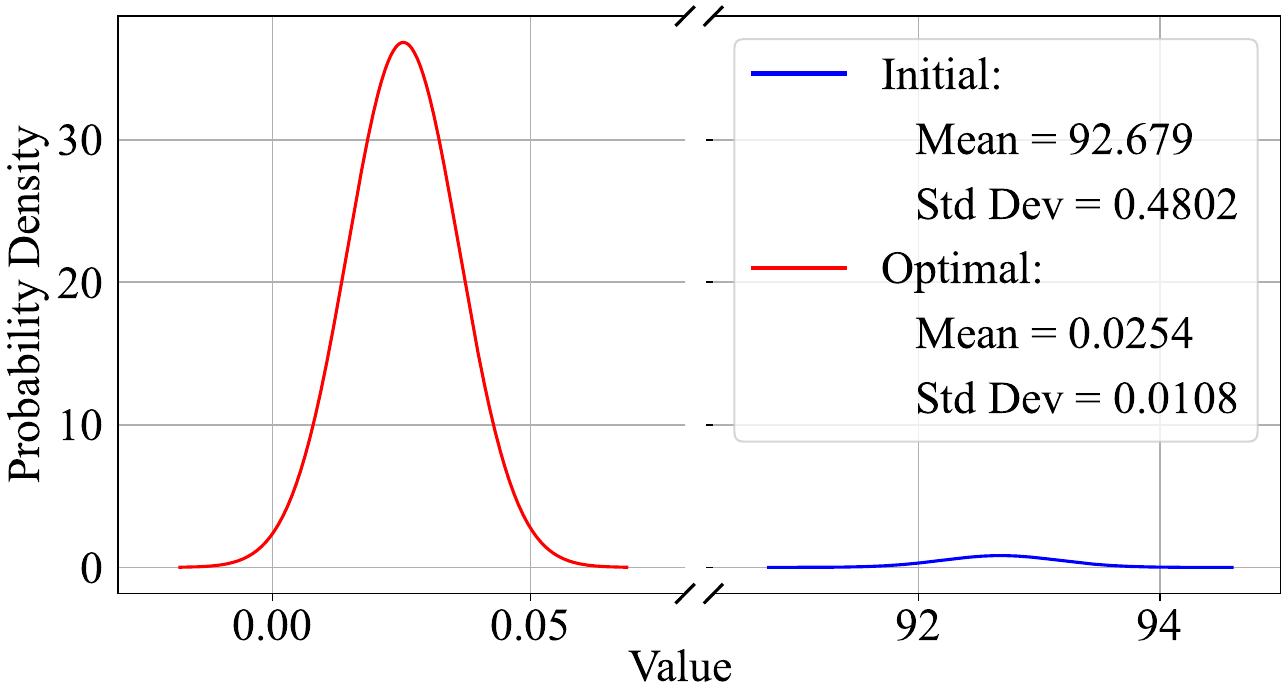}
        \caption*{(c) Case 2 (10D) - BNN-PDD MPSS}
    \end{minipage}
    \hfill
    \begin{minipage}[t]{0.48\textwidth}
        \centering
        \includegraphics[width=\textwidth]{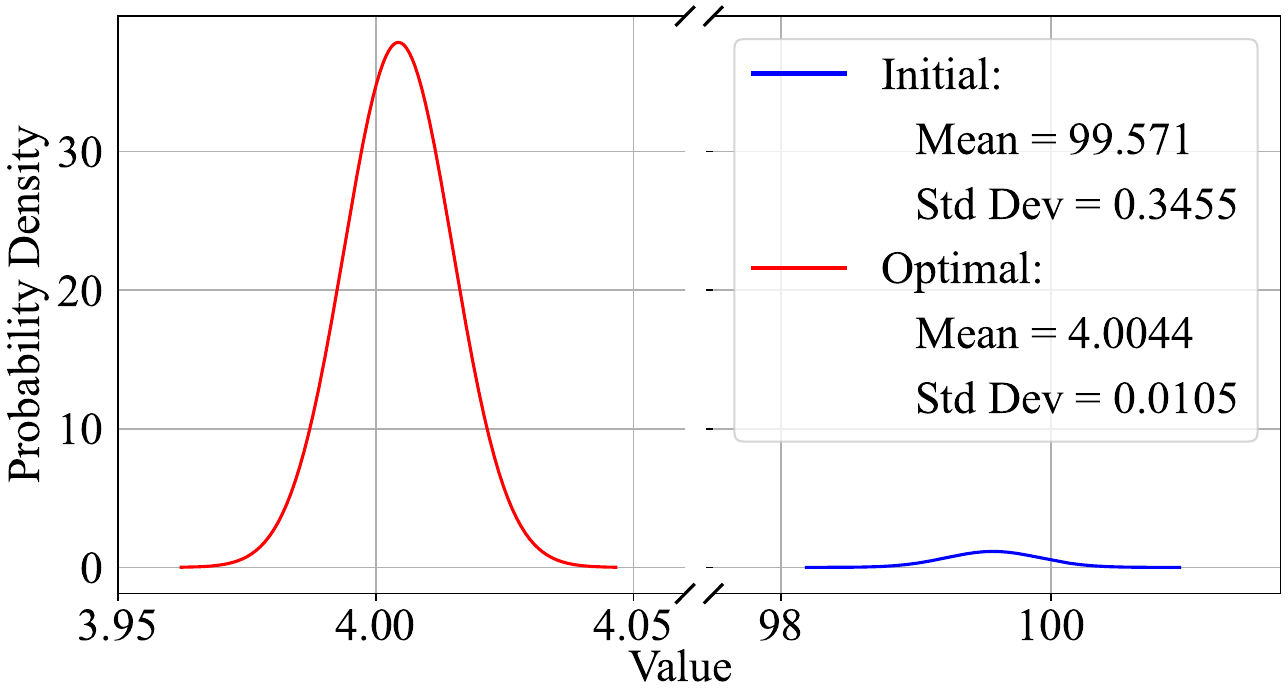}
        \caption*{(d) Case 2 (10D) - GP-PDD MPSS}
    \end{minipage}
    \caption{Comparison of PDF results for initial (blue) and optimal (red) designs between the proposed BNN-PDD MPSS and GP-PDD MPSS: (a) BNN-PDD MPSS results for Case 1 (2D), (b) GP-PDD MPSS results for Case 1 (2D), (c) BNN-PDD MPSS results for Case 2 (10D), (d) GP-PDD MPSS results for Case 2 (10D).}
    \label{fig:Rastrigin_pdf}
    \end{figure}

    Figure~\ref{fig:Rastrigin_pdf}(a) and (b) present the PDFs for the initial and optimal designs in Case~1. Both methods produce sharp peaks centered at zero in the optimal design distributions, in contrast to the initial design distributions. This result is attributed to the relatively high data density in the low-dimensional space, which provides sufficient prediction accuracy for both surrogate models. 
    Table~\ref{table:rastrigin_results} confirms that both methods achieve the global optimum, reducing the mean to approximately $0.0043$.

    \begin{table}[htb!]
    \centering
    \renewcommand{\arraystretch}{1.3}
    \caption{Optimization results for Example 1}
    \label{table:rastrigin_results}
    \begin{tabular}{@{}lllllll@{}}
    \toprule
     & Method & Statistics & Initial ($\mathbf{d}_{0}^{(0)}$) & Optimal ($\mathbf{d}^*$) & Reduction ($\%$) & No. of $y$ evaluations \\ \midrule
    Case 1 & BNN-PDD MPSS & Mean & 10.3782 & 0.0043 & \, 99.96 & 8300 \\
    (2D) &  & Std Dev & 2.5393$\times 10^{-1}$ & 4.1584$\times 10^{-3}$ & \, 98.36 &  \\ \cmidrule(l){2-7} 
     & GP-PDD MPSS & Mean & 5.6669 & 0.0043 & \, 99.92 & 4800 \\
     &  & Std Dev & 1.6034$\times 10^{-1}$ & 4.1399$\times 10^{-3}$ & \, 97.42 &  \\ \cmidrule(l){2-7} 
     & direct MCS & Mean & \, - & 0.03884 & \, - & 75000 \\
     &  & Std Dev & \, - & 5.1674$\times 10^{-2}$ & \, - &  \\ \midrule
    Case 2 & BNN-PDD MPSS & Mean & 92.6790 & 0.0254 & \, 99.97 & 16400 \\
    (10D) &  & Std Dev & 4.8018$\times 10^{-1}$ & 1.0831$\times 10^{-2}$ & \, 97.78 &  \\ \cmidrule(l){2-7} 
     & GP-PDD MPSS & Mean & 99.571 & 4.0044 & \, 95.98 & 20300 \\
     &  & Std Dev & 3.4549$\times 10^{-1}$ & 1.0536$\times 10^{-2}$ & \, 96.95 &  \\ \cmidrule(l){2-7} 
     & direct MCS & Mean & \, - & 16.1134 & \, - & 1.977$\times 10^{6}$ \\
     &  & Std Dev & \, - & 2.5816$\times 10^{-1}$ & \, - &  \\ \bottomrule
    \end{tabular}
    \end{table}
%

%%%%%%%%%%%%%%%%%%%%%%%%%%%%%%%%%%%%%%%%%%%
%
\subsubsection{Results for Case 2 ($N=M=10$)} \label{example 1_result_case2}
    Unlike Case 1, the high-dimensional problem in Case 2 reveals clear performance differences among the three methods. Although Fig~\ref{fig:Rastrigin_obj}(b) shows that both BNN-PDD MPSS and GP-PDD MPSS reduce the objective function value over iterations, Fig~\ref{fig:Rastrigin_pdf} shows a distinct difference. As shown in Fig~\ref{fig:Rastrigin_pdf}(d), the GP-based method reduces the standard deviation by $96.95\%$, yet the mean remained at $4.0044$, producing a distribution biased away from zero. This indicates that the GP kernel performance degrades with increasing dimensionality, reducing prediction accuracy and causing the optimizer to become trapped in a local optimum. In contrast, Fig~\ref{fig:Rastrigin_pdf}(c) demonstrates the robustness of the BNN-based method, where the PDF remains centered at zero even in the high-dimensional setting. This confirms that BNN effectively captures complex variable interactions, enabling accurate exploration of the global optimum in high-dimensional spaces. MCS, meanwhile, fails to identify a meaningful optimal solution within the allowed computational budget. The cure of dimensionality causes the search space to expand exponentially, making simple random search unable to locate the global optimum of the highly multimodal Rastrigin function.
    
    Table~\ref{table:rastrigin_results} summarizes the quantitative comparison. In Case 2, the GP-based method yields a final mean of $4.0044$, while MCS produces a mean of 16.1134, both failing to reach the global optimum. In contrast, BNN-PDD MPSS reduces the mean to $0.0254$ and the standard deviation to $1.0831 \times 10^{-2}$, confirming its ability to accurately and robustly locate the global optimum even in high-dimensional problems.

%
%%%%%%%%%%%%%%%%%%%%%%%%%%%%%%%%%%%%%%%%%%%
%
\subsection{Example 2: Shape design of the electric motor} \label{section:SPMSM}
    The second example addresses the shape optimization of an external rotor permanent magnet synchronous motor (ERPMSM). In this configuration, the rotor is located outside the stator, which secures a larger air gap radius compared to conventional inner rotor motors and thus delivers higher torque density within the same volume. This motor type is widely adopted in applications demanding high torque with compact form factors and constant-speed operation, such as drones, cooling fans, and robotics.
    
   Cogging torque arises from the magnetic attraction between the permanent magnets and the stator core even in the absence of excitation, and constitutes the primary source of motor vibration and noise. Minimizing cogging torque through design optimization is therefore critical. However, the complex nonlinear relationships between design variables and motor performance require FEA, which imposes high computational cost and makes RDO particularly challenging. Building on the validation results from Section~\ref{section:Rastrigin function}, this example evaluates the efficiency and robustness of the proposed BNN-PDD MPSS method as a practical motor design optimization strategy.
%

%
%%%%%%%%%%%%%%%%%%%%%%%%%%%%%%%%%%%%%%%%%%%
%
\subsubsection{Problem definition} \label{section:example2_problem definition}
  \begin{figure}[htb!]
        \centering
        \includegraphics[width=0.5\textwidth, trim=0 330 100 70]{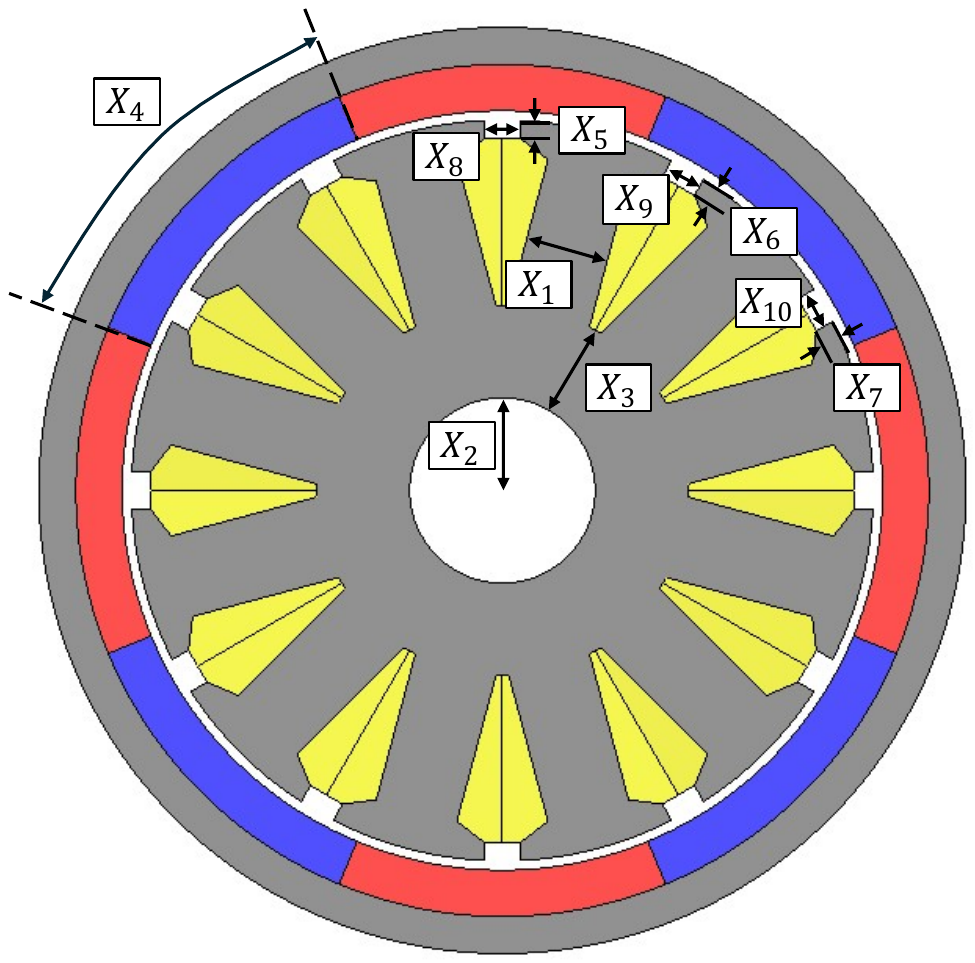}
        \caption{2D electromagnetic field analysis model of 8-pole 12-slot ERPMSM with random variables. Red and blue represent N and S poles of permanent magnets, yellow indicates the coil winding area, and gray represents iron core.}
        \label{fig:ERPMSM 2D model}
    \end{figure}
    The ERPMSM model used in this example is a surface-mounted permanent magnet synchronous motor with 8 poles and 12 slots. The poles refer to the number of permanent magnets attached to the rotor that generate magnetic flux, while the slots denote the number of grooves accommodating the stator windings. Figure~\ref{fig:ERPMSM 2D model} illustrates the two-dimensional electromagnetic field analysis model. We verified the FEA model through a comparative study  with~\cite{ling2016influence}: a cogging torque comparison at a magnet pole arc corresponding to a full pole pitch yielded a relative error of approximately $0.69 \%$.
    
    The motor operates at the rated speed of 750 rpm with an axial length of 50 mm. We adopt neodymium magnets with a relative permeability of 1.05 and a remanence of $0.85 \mathrm{T}$ and apply AISI 1010 steel for the iron core.
    
    We leverage the geometric symmetry of the motor to efficiently model manufacturing tolerances. In the 8-pole, 12-slot topology, the greatest common divisor of the pole and slot numbers is four, indicating that both the electromagnetic flux distribution and the geometry repeat four times per mechanical revolution. We therefore construct a one-quarter model, which represents the minimal periodic unit. We define random variables exclusively for the geometric dimensions of the first quarter sector. The remaining three sectors periodically replicate these values.
    
    To incorporate manufacturing uncertainties, we model the geometric dimensions as random variables. We define a total of 10 random variables $\mathbf{X} = (X_1, \ldots, X_{10})^{\intercal}$, each following a truncated Gaussian distribution. For the stator tooth dimensions $d_5$ and $d_6$, independent manufacturing tolerances can occur at each of the three slots within the one-quarter model. We therefore model these dimensions using three independent random variables each $X_5 \sim X_7$ and $X_8 \sim X_{10}$, respectively. Table~\ref{table:X} summarizes the random variables, manufacturing tolerances, and design ranges. We define the design variables as the mean values of the corresponding random variable. The six design variables are $d_k = \mathbb{E}[\mathbf{X}_{k}]$ for $k = 1,2,3,4$, $d_5 = \mathbb{E}[\mathbf{X}_5] = \mathbb{E}[\mathbf{X}_6] = \mathbb{E}[\mathbf{X}_7]$, and $d_6 = \mathbb{E}[\mathbf{X}_8] = \mathbb{E}[\mathbf{X}_9] = \mathbb{E}[\mathbf{X}_{10}]$.
    \begin{table}[htb!]
    \caption{Statistical properties of random variables of ERPMSM (Example 2)}\label{table:X}
    \centering
    \renewcommand{\arraystretch}{1.2}
    \begin{tabular}{@{}llllll@{}}
    \toprule
    Random variables & \begin{tabular}[c]{@{}l@{}}Design variables\\ (Mean)\end{tabular} & Tolerance & Min. & Max. & \begin{tabular}[c]{@{}l@{}}Probability\\ distribution\end{tabular} \\ \midrule
    \begin{tabular}[c]{@{}l@{}}Tooth body thickness\\ ($X_1$), mm\end{tabular} & \, $d_1$ & \, 0.04 & \, 2.5 & \, 5.0 & Truncated Gaussian \\
    \begin{tabular}[c]{@{}l@{}}Stator inner radius\\ ($X_2$), mm\end{tabular} & \, $d_2$ & \, 0.04 & \, 5.0 & \, 7.5 & Truncated Gaussian \\
    \begin{tabular}[c]{@{}l@{}}Stator back yoke thickness\\ ($X_3$), mm\end{tabular} & \, $d_3$ & \, 0.04 & \, 3.0 & \, 5.0 & Truncated Gaussian \\
    \begin{tabular}[c]{@{}l@{}}PM pole arc\\ ($X_4$), $^\circ$ \end{tabular} & \, $d_4$ & \, 0.04 & \, 31.0 & \, 45.0 & Truncated Gaussian \\
    \begin{tabular}[c]{@{}l@{}}Tooth tip thickness \\ ($X_5, X_6, X_7$), mm\end{tabular} & \, $d_5$ & \, 0.04 & \, 0.3 & \, 1.0 & Truncated Gaussian \\
    \begin{tabular}[c]{@{}l@{}}Slot opening width\\ ($X_8, X_9, X_{10}$), mm\end{tabular} & \, $d_6$ & \, 0.04 & \, 0.5 & \, 3.0 & Truncated Gaussian \\ \bottomrule
    \end{tabular}
    \end{table}
        
    We evaluate cogging torque through electromagnetic field analysis under a no-load condition as cogging torque arises without electrical excitation. We formulate the RDO problem to simultaneously minimize the mean and standard deviation of the cogging torque. We impose a constraint to ensure that the mean of cogging torque does not exceed the allowable limit $T_{cogging}=0.39 \; \mathrm{Nm}$. This RDO is defined as  
    \begin{equation*}
        \begin{aligned}
            \min\limits_{\mathbf{d} \in \mathcal{D}} \quad &c_{0} = w_{1} \frac{\mathbb{E}_{\mathbf{d}}[y(\mathbf{X})]}{\mathbb{E}_{\mathbf{d}_{0}}[y(\mathbf{X})]} + w_{2} \frac{\sqrt{\mathbb{V}\mathrm{ar}_{\mathbf{d}}[y(\mathbf{X})]}}{\sqrt{\mathbb{V}\mathrm{ar}_{\mathbf{d}_{0}}[y(\mathbf{X})]}}, \\
            \text{subject to} \quad &c_{1} = \frac{\mathbb{E}_{\mathbf{d}}[y(\mathbf{X})] - T_{cogging}}{\mathbb{E}_{\mathbf{d}_{0}}[y(^{*}\mathbf{X;d})]} \leq 0, \\
            &2.5 \, \text{mm} \leq d_1 \leq 5.0 \, \text{mm}, \\
            &5.0 \, \text{mm} \leq d_2 \leq 7.5 \, \text{mm}, \\
            &3.0 \, \text{mm} \leq d_3 \leq 5.0 \, \text{mm}, \\
            &31.0 \, ^\circ \leq d_4 \leq 45.0 \, ^\circ, \\
            &0.3 \, \text{mm} \leq d_5 \leq 1.0 \, \text{mm}, \\
            &0.5 \, \text{mm} \leq d_6 \leq 3.0 \, \text{mm}, \\
        \end{aligned}
    \end{equation*}
    where $\mathbb{E}_{\mathbf{d}_{0}}[y(\mathbf{X})]$ and $\mathbb{V}\mathrm{ar}_{\mathbf{d}_{0}}[y(\mathbf{X})]$ denote the mean and variance of cogging torque at the initial design $\mathbf{d}_{0}$. The initial design vector is $\mathbf{d}_{0} = (4.5 \, \mathrm{mm}, 5.0 \, \mathrm{mm}, 5.0 \, \mathrm{mm}, 45.0^\circ, 1.0 \, \mathrm{mm}, 2.0 \, \mathrm{mm})^{\intercal}$, and the weights are set to $w_1 = w_2 = 0.5$.

\subsubsection{Method} \label{example2_method}
    We apply the proposed BNN-PDD MPSS to construct local surrogate models in each subregion. We generate 500 training samples within the six-dimensional design space for each subregion using LHS. 
%
%%%%%%%%%%%%%%%%%%%%%%%%%%%%%%%%%%%%%%%%%%%
%
\subsubsection{Results} \label{section:example2_results}
    \begin{figure}[htb!]
    \centering
    \begin{minipage}[t]{0.48\textwidth}
        \centering
        \includegraphics[width=\textwidth]{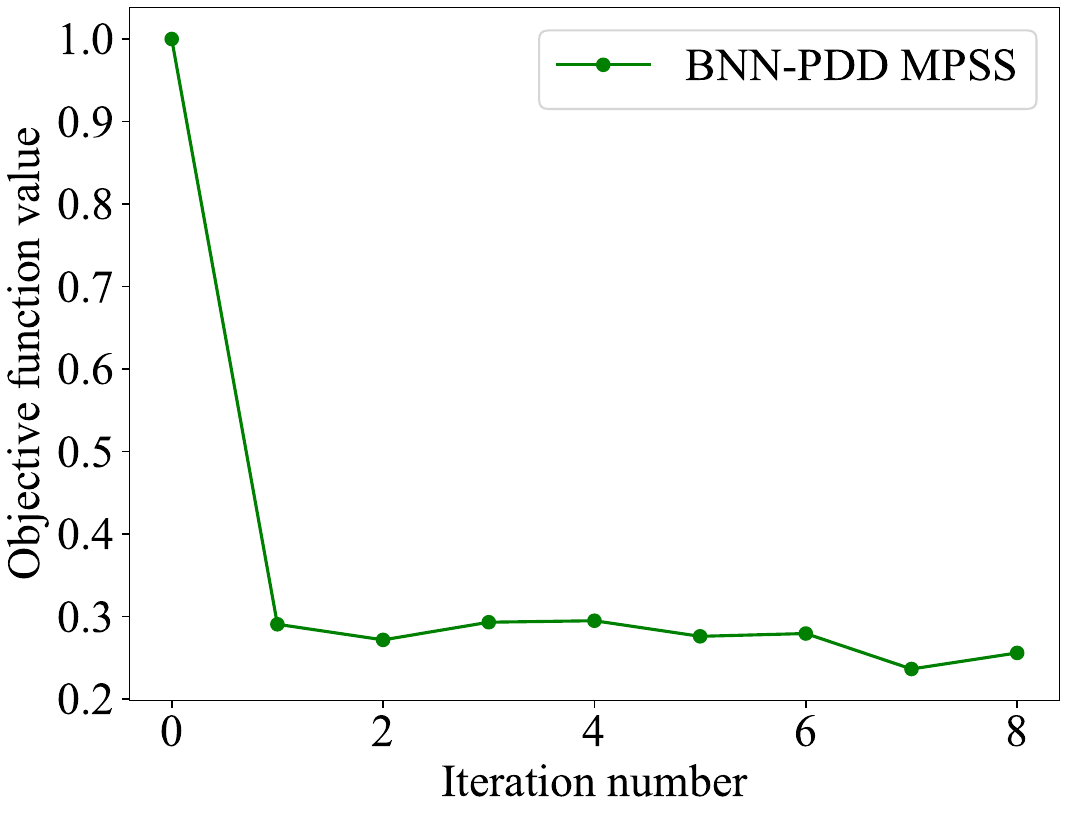}
        \caption*{(a)}
    \end{minipage}
    \hfill
    \begin{minipage}[t]{0.48\textwidth}
        \centering
        \includegraphics[width=\textwidth]{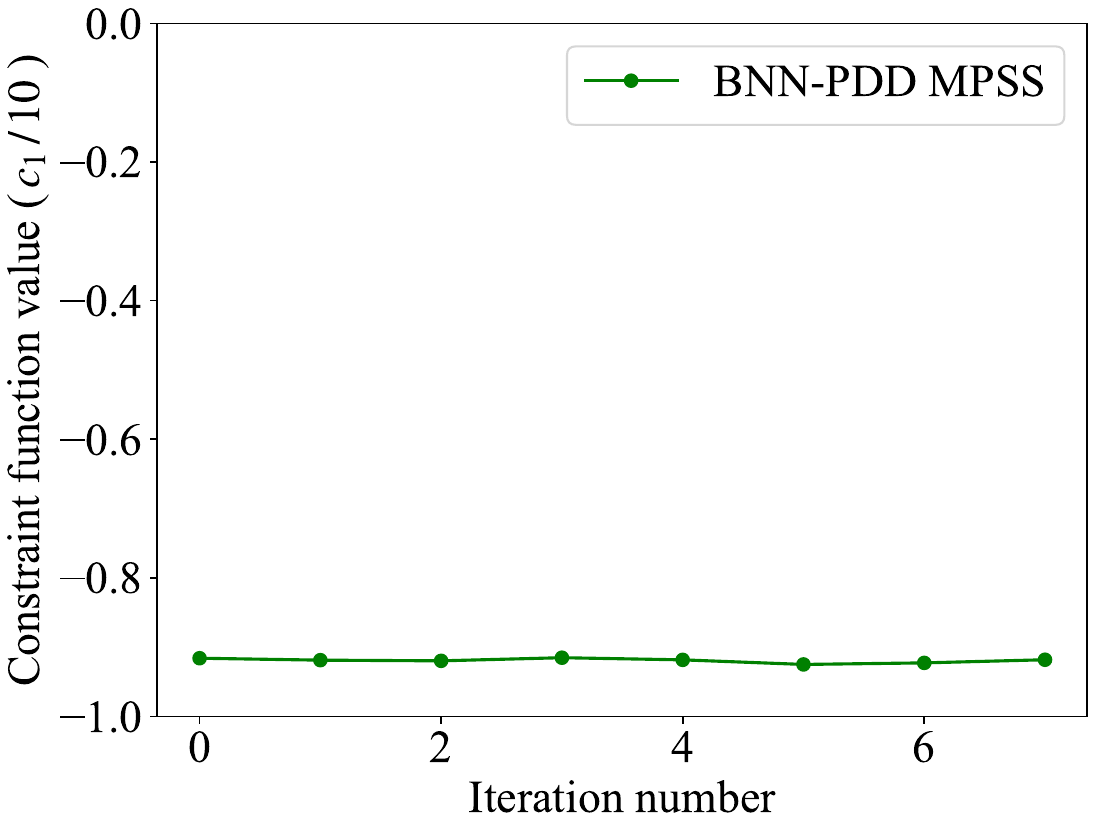}
        \caption*{(b)}
    \end{minipage}
    \caption{RDO iteration history for Example 2: (a) objective function, and (b) constraint function}
    \label{fig:ERPMSM_BNN_obj}
    \end{figure}
    \begin{figure}[htb!]
        \centering
        \includegraphics[width=0.8\textwidth]{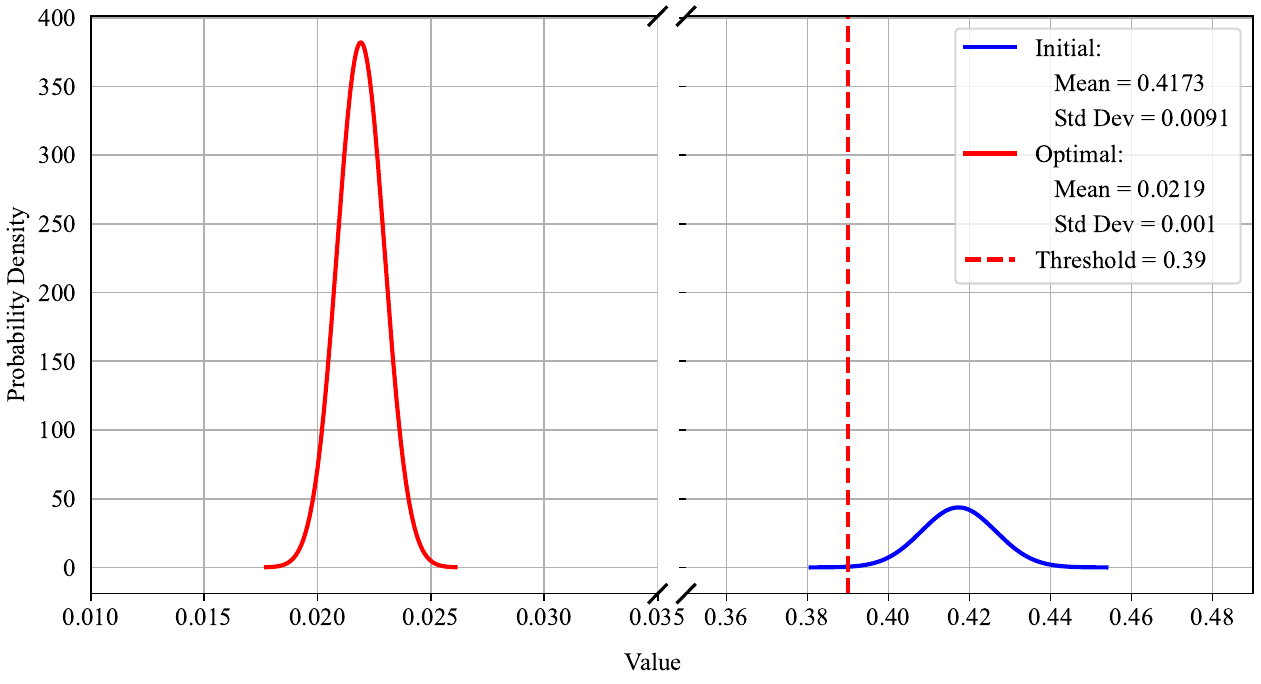}
        \caption{Comparison of PDF results for initial (blue) and optimal (red) designs for Example 2}
        \label{fig:ERPMSM_BNN_pdf}
    \end{figure}
    
    Figure~\ref{fig:ERPMSM_BNN_obj}(a) and (b) show the iteration histories of the objective and constraint functions across subregion iterations. The objective function decreases monotonically, indicating stable convergence. Figure~\ref{fig:ERPMSM_BNN_pdf} compares the PDFs of the cogging torque for the initial and optimal designs. The initial design (blue) exhibits a broad distribution centered around $0.42$~Nm, lying entirely above the allowable threshold $T_{\text{cog}} = 0.39$~Nm (red dashed line) and thus violating the constraint. In contrast, the optimal design (red) produces a sharp, narrow peak concentrated near $0.022$~Nm---well below the threshold---confirming that the proposed method simultaneously reduces the mean cogging torque and suppresses its variability while satisfying the constraint. 
    \begin{table}[htb!]
    \caption{Optimization results for Example 2}
    \label{table:example2_results}
    \centering
    \renewcommand{\arraystretch}{1.2}
    \begin{tabular}{@{}lllll@{}}
    \toprule
    Statistics & Initial ($\mathbf{d}_{0}$) & Optimal ($\mathbf{d}^*$) & Reduction ($\%$) & No. of FEA \\ \midrule
    Mean {[}Nm{]} & 0.4173 & 0.02191 & \, 94.75 & 6644 \\
    Std Dev {[}Nm{]} & 9.1274$\times 10^{-3}$ & 1.0444$\times 10^{-3}$ & \, 88.56 &  \\ \bottomrule
    \end{tabular}
    \end{table}
     \begin{figure}[htb!]
    \centering
    \begin{minipage}[t]{0.43\textwidth}
        \centering
        \includegraphics[width=\textwidth]{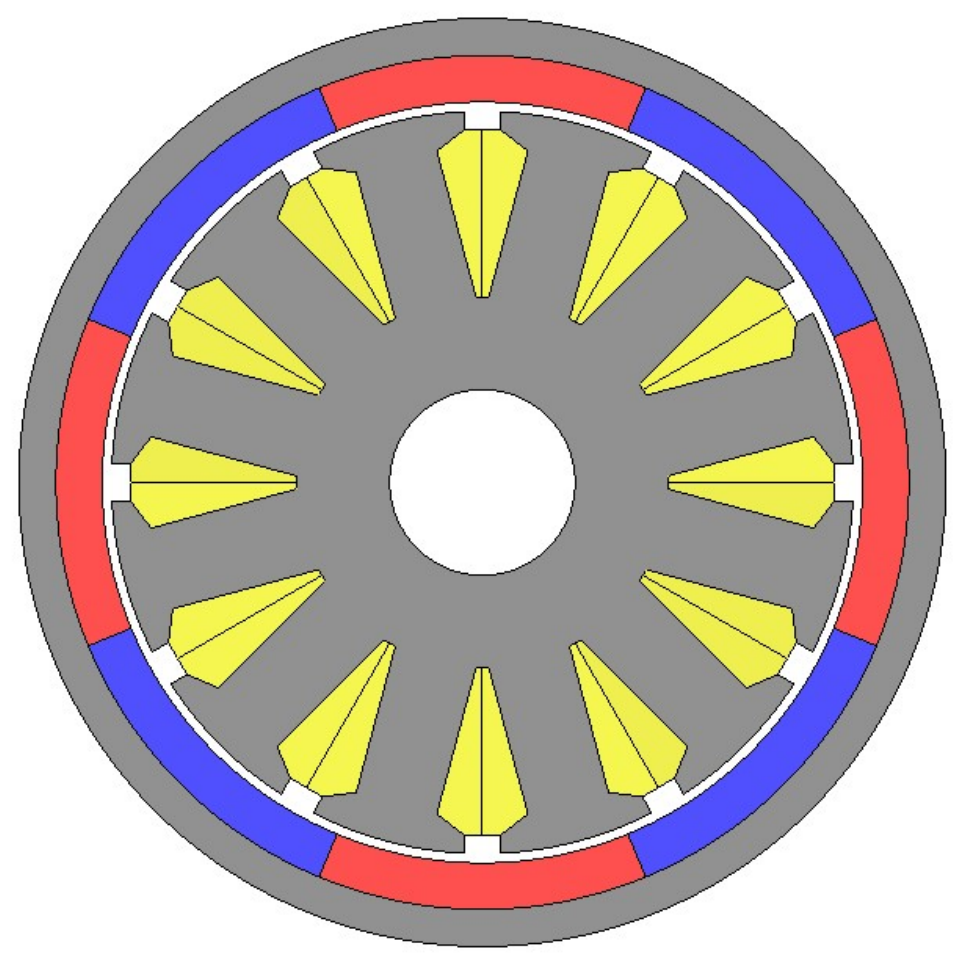}
        \caption*{(a)}
    \end{minipage}
    \hfill
    \begin{minipage}[t]{0.45\textwidth}
        \centering
        \includegraphics[width=\textwidth]{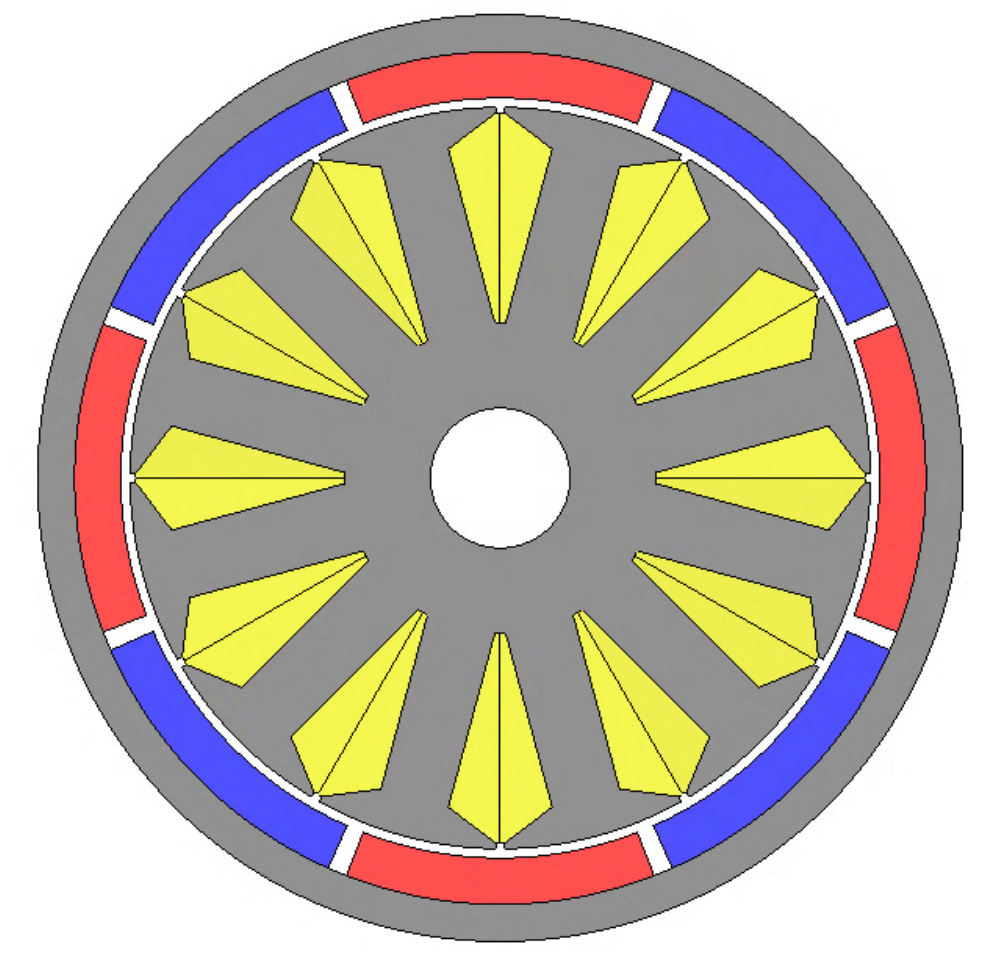}
        \caption*{(b)}
    \end{minipage}
    \caption{Comparison of the 2D topology between (a) initial design and (b) optimal design of the ERPMSM.}
    \label{fig:ERPMSM_model}
    \end{figure}   
    Table~\ref{table:example2_results} summarizes the quantitative results. The proposed BNN-PDD MPSS method reduces the mean cogging torque from $0.4173$~Nm to $0.02191$~Nm ($94.75\%$ reduction) and the standard deviation from $9.1274 \times 10^{-3}$~Nm to $1.0444 \times 10^{-3}$~Nm ($88.56\%$ reduction), achieving both performance improvement and robustness enhancement with only 6644 FEA evaluations. This result confirms that the proposed method effectively handles complex electromagnetic field problems, delivering a robust optimal design even in a high-dimensional, strongly nonlinear setting. Figure~\ref{fig:ERPMSM_model} compares the 2D topologies of the initial and optimal designs.
    %Table~\ref{table:example2_results} summarizes the quantitative results, including statistical moments and reduction rates. The initial design shows a mean of $0.4173 \, \mathrm{Nm}$ and a standard deviation of $9.1274 \times 10^{-3} \, \mathrm{Nm}$. In contrast, the optimal design achieves significant improvement, reaching a mean of $0.02191 \, \mathrm{Nm}$ and a standard deviation of $1.0444 \times 10^{-3} \, \mathrm{Nm}$. Consequently, the optimal design reduces the mean by $94.75 \%$ and the standard deviation by $88.56 \%$ compared to the initial design. This result demonstrates that the proposed BNN-based RDO methodology secures both the performance and robustness of cogging torque, even in complex electromagnetic field problems. Figure~\ref{fig:ERPMSM_model} shows the 2D topologies of the initial and optimal designs.

%

%
%%%%%%%%%%%%%%%%%%%%%%%%%%%%%%%%%%%%%%%%%%%
\section{Conclusion} \label{section:conclusion}
%%%%%%%%%%%%%%%%%%%%%%%%%%%%%%%%%%%%%%%%%%%
%
   This study proposed a novel RDO method that integrates BNN with PDD. The combined method leverages BNN as a surrogate to replace expensive model evaluations required for PDD coefficient calculation, enabling rapid estimation of statistical moments at significantly reduced computational cost. The proposed method maximizes sampling efficiency through an active learning strategy driven by BNN predictive uncertainty. Adopting a MPSS strategy, the method partitions the design space into subregions and constructs local surrogate models, while employing a global search within each subregion to handle severe nonlinearity. A dynamic resizing strategy further adjusts subregion boundaries during optimization to avoid entrapment in local optima and enhance search robustness. PDD is then applied to the resulting surrogate, exploiting the structural decomposition of component functions for efficient and accurate probabilistic analysis.

Numerical results confirmed that the proposed method converges stably to optimal solutions with significantly fewer function evaluations compared to conventional approaches. In the high-dimensional Rastrigin problem ($N=10$), BNN-PDD MPSS achieved a $99.97\%$ mean reduction, while GP-PDD MPSS and MCS failed to locate the global optimum. In the ERPMSM design problem, the method reduced the mean cogging torque by $94.75\%$ and the standard deviation by $88.56\%$ with only 6644 FEA evaluations, demonstrating its practicality for high-dimensional, strongly nonlinear engineering problems.

The proposed method provides an effective solution for robust design of complex engineering systems requiring high-fidelity analysis. However, BNN has inherent limitations: the training process is more complex than that of conventional neural networks, and determining optimal hyperparameters and prior distributions involves relatively high computational cost and sensitivity. The Bayesian optimization-based tuning used in this study required iterative model retraining to identify optimal hyperparameters, resulting in substantial computational overhead. Future work will focus on introducing algorithms that simultaneously optimize weights and hyperparameters during the training process, eliminating the need for iterative retraining and further improving the overall efficiency of the method.

\begin{appendices}
%%%%%%%%%%%%%%%%%%%%%%%%%%%%%%%%%%%%%%%%%%%

\printnomenclature

\end{appendices}
%%%%%%%%%%%%%%%%%%%%%%%%%%%%%%%%%%%%%%%
%%%%%%%%%%%%%%%%%%%%%%%%%%%%%%%%%%%%%%%%
% BIBLIOGRAPHY
%%%%%%%%%%%%%%%%%%%%%%%%%%%%%%%%%%%%%%%

\section*{Acknowledgment} 
This work was supported by Institute of Information \& Communications Technology Planning \& Evaluation (IITP) grant funded by the Korea government (MSIP) (No. RS-2025-02304285, Development of Digital Triplet Based Predictive Maintenance Solution for Future Power Facility Innovation).

\section*{Conflict of Interest}
On behalf of all authors, the corresponding author states that there is no conflict of interest.

\section*{Replication of Results}
All BNN architecture details, PDD truncation parameters, MPSS tolerances, and active learning criteria are fully specified in Section~\ref{section:algorithm}. Due to the proprietary nature of the FEA model used in the motor design example, the simulation files cannot be made publicly available. The source code for the proposed BNN-PDD MPSS method is available from the corresponding author upon request.
\bibliography{references}

@book{taguchi_taguchi_1993,
	title = {Taguchi on {Robust} {Technology} {Development}: {Bringing} {Quality} {Engineering} {Upstream}},
	isbn = {0-7918-0028-8},
	url = {https://doi.org/10.1115/1.800288},
	publisher = {ASME Press},
	author = {Taguchi, Genichi},
	month = jan,
	year = {1993},
	doi = {10.1115/1.800288},
}

@article{mourelatos2006methodology,
  title={A methodology for trading-off performance and robustness under uncertainty},
  author={Mourelatos, Zissimos P and Liang, Jinghong},
  journal={J. Mech. Des.},
  volume={128},
  number={4},
  pages={856--863},
  year={2006}
}

@article{zaman2011robustness,
  title={Robustness-based design optimization under data uncertainty},
  author={Zaman, Kais and McDonald, Mark and Mahadevan, Sankaran and Green, Lawrence},
  journal={Structural and Multidisciplinary Optimization},
  volume={44},
  number={2},
  pages={183--197},
  year={2011},
  publisher={Springer}
}

@article{park2006robust,
  title={Robust design: an overview},
  author={Park, Gyung-Jin and Lee, Tae-Hee and Lee, Kwon Hee and Hwang, Kwang-Hyeon},
  journal={AIAA journal},
  volume={44},
  number={1},
  pages={181--191},
  year={2006}
}

@article{yao2011review,
  title={Review of uncertainty-based multidisciplinary design optimization methods for aerospace vehicles},
  author={Yao, Wen and Chen, Xiaoqian and Luo, Wencai and Van Tooren, Michel and Guo, Jian},
  journal={Progress in Aerospace Sciences},
  volume={47},
  number={6},
  pages={450--479},
  year={2011},
  publisher={Elsevier}
}

@article{kim2025interpolation,
  title={Interpolation-based optimal knot selection in spline dimensional decomposition for uncertainty quantification in dynamical systems},
  author={Kim, Yeonsu and Lee, Junhan and Wang, Bingran and Hwang, John T and Lee, Dongjin},
  journal={Applied Mathematical Modelling},
  pages={116613},
  year={2025},
  publisher={Elsevier}
}

@article{lee2023multifidelity,
  title={Multifidelity conditional value-at-risk estimation by dimensionally decomposed generalized polynomial chaos-{K}riging},
  author={Lee, Dongjin and Kramer, Boris},
  journal={Reliability Engineering \& System Safety},
  volume={235},
  pages={109208},
  year={2023},
  publisher={Elsevier}
}

@article{guo2025robust,
  title={Robust design optimization with limited data for char combustion},
  author={Guo, Yulin and Lee, Dongjin and Kramer, Boris},
  journal={Structural and Multidisciplinary Optimization},
  volume={68},
  number={3},
  pages={1--18},
  year={2025},
  publisher={Springer}
}

@article{hauzenberger2025bayesian,
  title={Bayesian neural networks for macroeconomic analysis},
  author={Hauzenberger, Niko and Huber, Florian and Klieber, Karin and Marcellino, Massimiliano},
  journal={Journal of Econometrics},
  volume={249},
  pages={105843},
  year={2025},
  publisher={Elsevier}
}

@article{song2025ferroelectric,
  title={Ferroelectric NAND for efficient hardware bayesian neural networks},
  author={Song, Minsuk and Koo, Ryun-Han and Kim, Jangsaeng and Han, Chang-Hyeon and Yim, Jiyong and Ko, Jonghyun and Yoo, Sijung and Choe, Duk-hyun and Kim, Sangwook and Shin, Wonjun and others},
  journal={Nature Communications},
  volume={16},
  number={1},
  pages={6879},
  year={2025},
  publisher={Nature Publishing Group UK London}
}

@article{ngartera2024application,
  title={Application of bayesian neural networks in healthcare: three case studies},
  author={Ngartera, Lebede and Issaka, Mahamat Ali and Nadarajah, Saralees},
  journal={Machine Learning and Knowledge Extraction},
  volume={6},
  number={4},
  pages={2639--2658},
  year={2024},
  publisher={MDPI}
}

@incollection{goan2020bayesian,
  title={Bayesian neural networks: An introduction and survey},
  author={Goan, Ethan and Fookes, Clinton},
  booktitle={Case Studies in Applied Bayesian Data Science: CIRM Jean-Morlet Chair, Fall 2018},
  pages={45--87},
  year={2020},
  publisher={Springer}
}

@article{mackay1992practical,
  title={A practical Bayesian framework for backpropagation networks},
  author={MacKay, David JC},
  journal={Neural computation},
  volume={4},
  number={3},
  pages={448--472},
  year={1992},
  publisher={MIT Press One Rogers Street, Cambridge, MA 02142-1209, USA journals-info~…}
}

@article{rahman2008polynomial,
  title={A polynomial dimensional decomposition for stochastic computing},
  author={Rahman, Sharif},
  journal={International Journal for Numerical Methods in Engineering},
  volume={76},
  number={13},
  pages={2091--2116},
  year={2008},
  publisher={Wiley Online Library}
}

@article{rahman2010statistical,
  title={Statistical moments of polynomial dimensional decomposition},
  author={Rahman, Sharif},
  journal={Journal of engineering mechanics},
  volume={136},
  number={7},
  pages={923--927},
  year={2010},
  publisher={American Society of Civil Engineers}
}

@article{jordan1999introduction,
  title={An introduction to variational methods for graphical models},
  author={Jordan, Michael I and Ghahramani, Zoubin and Jaakkola, Tommi S and Saul, Lawrence K},
  journal={Machine learning},
  volume={37},
  number={2},
  pages={183--233},
  year={1999},
  publisher={Springer}
}

@article{graves2011practical,
  title={Practical variational inference for neural networks},
  author={Graves, Alex},
  journal={Advances in neural information processing systems},
  volume={24},
  year={2011}
}

@article{helton2003latin,
  title={Latin hypercube sampling and the propagation of uncertainty in analyses of complex systems},
  author={Helton, Jon C and Davis, Freddie Joe},
  journal={Reliability Engineering \& System Safety},
  volume={81},
  number={1},
  pages={23--69},
  year={2003},
  publisher={Elsevier}
}

@article{janssen2013monte,
  title={Monte-Carlo based uncertainty analysis: Sampling efficiency and sampling convergence},
  author={Janssen, Hans},
  journal={Reliability Engineering \& System Safety},
  volume={109},
  pages={123--132},
  year={2013},
  publisher={Elsevier}
}

@article{toropov1993multiparameter,
  title={Multiparameter structural optimization using FEM and multipoint explicit approximations},
  author={Toropov, Vassili V and Filatov, AA and Polynkin, AA},
  journal={Structural optimization},
  volume={6},
  number={1},
  pages={7--14},
  year={1993},
  publisher={Springer}
}

@inproceedings{ling2016influence,
  title={Influence of magnet pole arc variation on the performance of external rotor permanent magnet synchronous machine based on finite element analysis},
  author={Ling, PP and Ishak, Dahaman and Tiang, TL},
  booktitle={2016 IEEE International Conference on Power and Energy (PECon)},
  pages={552--557},
  year={2016},
  organization={IEEE}
}

@article{lee2021robust,
  title={Robust design optimization under dependent random variables by a generalized polynomial chaos expansion},
  author={Lee, Dongjin and Rahman, Sharif},
  journal={Structural and Multidisciplinary Optimization},
  volume={63},
  number={5},
  pages={2425--2457},
  year={2021},
  publisher={Springer}
}

@article{guo2025risk,
  title={Risk-based design optimization for powder bed fusion metal additive manufacturing},
  author={Guo, Yulin and Kramer, Boris},
  journal={Structural and Multidisciplinary Optimization},
  volume={68},
  number={9},
  pages={179},
  year={2025},
  publisher={Springer}
}

@article{lee2020practical,
  title={Practical uncertainty quantification analysis involving statistically dependent random variables},
  author={Lee, Dongjin and Rahman, Sharif},
  journal={Applied Mathematical Modelling},
  volume={84},
  pages={324--356},
  year={2020},
  publisher={Elsevier}
}

@article{cho2023performance,
  title={Performance estimation of freeze protection system for outdoor fire piping by using {AI} algorithm},
  author={Cho, Hojoon and Seo, Sangmin and Heo, Chinseok and Kwak, Junjae and Kim, Yongbae and Park, Jinsoo and Lee, Sangjun and Lim, Seongsik},
  journal={Journal of Mechanical Science and Technology},
  volume={37},
  number={10},
  pages={5093--5101},
  year={2023},
  publisher={Springer}
}

@article{hyun2026robust,
  title={Robust topology optimization of continuum structures under material uncertainties using multifidelity {M}onte {C}arlo},
  author={Hyun, Jaeyub and Chaudhuri, Anirban and Willcox, Karen E and Kim, H Alicia},
  journal={Computer Methods in Applied Mechanics and Engineering},
  volume={451},
  pages={118714},
  year={2026},
  publisher={Elsevier}
}

@article{lee2024novel,
  title={A novel sampling method for adaptive gradient-enhanced {K}riging},
  author={Lee, Mingyu and Noh, Yoojeong and Lee, Ikjin},
  journal={Computer Methods in Applied Mechanics and Engineering},
  volume={418},
  pages={116456},
  year={2024},
  publisher={Elsevier}
}

@article{kim4883196deep,
  title={Deep Reinforcement {N}elder-{M}ead Optimization for {H}vac Digital Twin Model Performance},
  author={Kim, Nuri and Noh, Yoojeong and Choi, Gyungmin and Oh, Dong Jun and Kang, Young-Jin and Park, Noma and Choi, Soonyong},
  journal={Available at SSRN 4883196}
}
\bibliographystyle{unsrt}

%\newpage
%\begin{appendix}
%    \input{appendix}
%\end{appendix}
%
\end{document}